\documentclass[11pt]{article}
\usepackage{amsmath} \usepackage{amssymb} \usepackage{latexsym} \usepackage{enumitem} \usepackage{mathrsfs} \usepackage{comment}
\usepackage{color}

\usepackage[colorlinks=true,urlcolor=blue,
citecolor=red,linkcolor=blue,linktocpage,pdfpagelabels,bookmarksnumbered,bookmarksopen]{hyperref}

\newtheorem{assumption}{Assumption}

\def\qed{ \ \vrule width.2cm height.2cm depth0cm\smallskip}
\newenvironment{proof}{\noindent {\bf Proof.\/}}{$\qed$\vskip 0.1in}

\newcommand{\la}{\langle}
\newcommand{\ra}{\rangle}

\newcommand{\ol}{\overline}
\newcommand{\ul}{\underline}

\newcommand{\ba}{\begin{array}}
\newcommand{\ea}{\end{array}}
\newcommand{\be}{\begin{equation}}
\newcommand{\ee}{\end{equation}}
\newcommand{\bea}{\begin{eqnarray}}
\newcommand{\eea}{\end{eqnarray}}
\newcommand{\beaa}{\begin{eqnarray*}}
\newcommand{\eeaa}{\end{eqnarray*}}

\def\dbE{\mathbb{E}}
\def\dbF{\mathbb{F}}
\def\dbG{\mathbb{G}}

\def\dbL{\mathbb{L}}

\def\dbP{\mathbb{P}}
\def\dbR{\mathbb{R}}

\def\sL{\mathscr{L}}

\def\a{\alpha}
\def\b{\beta}
\def\g{\gamma}
\def\d{\delta}
\def\e{\varepsilon}

\def\k{\kappa}
\def\l{\lambda}

\def\si{\sigma}

\def\f{\varphi}
\def\th{\theta}
\def\o{\omega}
\def\h{\widehat}
\def\G{\Gamma}
\def\D{\Delta}

\def\L{\Lambda}

\def\O{\Omega}
\def\cA{{\cal A}}

\def\cC{{\cal C}}

\def\cF{{\cal F}}

\def\cI{{\cal I}}

\def\cL{{\cal L}}
\def\cM{{\cal M}}

\def\cP{{\cal P}}

\def\cS{{\cal S}}

\def\cU{{\cal U}}

\def\cW{{\cal W}}

\def\no{\noindent}

\def\ss{\smallskip}
\def\ms{\medskip}

\def\q{\quad}
\def\qq{\qquad}

\def\pa{\partial}
\def\cd{\cdot}
\def\cds{\cdots}

\def\tr{\hbox{\rm tr}}

\def\qed{ \hfill \vrule width.25cm height.25cm depth0cm\smallskip}

\newcommand{\basa}{\begin{assumption}}
\newcommand{\easa}{\end{assumption}}

\newcommand{\bas}{\begin{assum}}
\newcommand{\eas}{\end{assum}}

\def\pa{\partial}
\def\h{\widehat}

 \def\cd{\cdot}
\def\cds{\cdots}

\def\tr{\hbox{\rm tr$\,$}}

\def\dis{\displaystyle}

\def\wh{\widehat}

\def\1{{\bf 1}}

\def\:{\!:\!}
\def\reff{\eqref}
\def \proof{{\noindent \bf Proof.\quad}}

at 9pt

\definecolor{alp}{rgb}{0.0, 0.5, 0.0}

\newtheorem{thm}{Theorem}[section]

\newtheorem{prop}[thm]{Proposition}
\newtheorem{rem}[thm]{Remark}
\newtheorem{eg}[thm]{Example}
\newtheorem{defn}[thm]{Definition}
\newtheorem{assum}[thm]{Assumption}

\begin{document}

\title{\bf Mean Field Games of Controls: 
Propagation of Monotonicities} 
\author{Chenchen Mou\thanks{\noindent  Dept. of Math.,
City University of Hong Kong. E-mail: \href{mailto:chencmou@cityu.edu.hk}{chencmou@cityu.edu.hk}. This author is supported in part by CityU Start-up Grant 7200684 and Hong Kong RGC Grant ECS 9048215.} ~ and ~ Jianfeng Zhang\thanks{\noindent  Dept. of Math., 
University of Southern California. E-mail:
\href{mailto:jianfenz@usc.edu}{jianfenz@usc.edu}. This author is supported in part by NSF grant DMS-1908665 and and DMS-2205972.
}  
}
\date{}
\maketitle

\begin{abstract} 
The theory of Mean Field Game of Controls considers a class of mean field games where the interaction is through  the joint distribution of the state and control. It is well known that, for standard mean field games, certain monotonicity condition is crucial to guarantee the uniqueness of mean field equilibria and then the global wellposedness for master equations. In the literature, the monotonicity condition could be the Lasry-Lions monotonicity, the displacement monotonicity, or the anti-monotonicity conditions. In this paper, we investigate all these three types of monotonicity conditions for Mean Field Games of Controls and show their propagation along the solutions to the master equations with common noises. In particular, we extend the displacement monotonicity to semi-monotonicity, whose propagation result\footnote{We would like to thank Gangbo and  M\'esz\'aros for helpful discussions on this part.} is new even for standard mean field games. This is the first step towards the global wellposedness theory for  master equations of Mean Field Games of Controls.
\end{abstract}

\no{\bf Keywords.}  Mean field game of controls, master equation, Lasry-Lions monotonicity, displacement semi-monotonicity, anti-monotonicity

\ms
\no{\it 2020 AMS Mathematics subject classification:}  35R15, 49N80, 60H30, 91A16, 93E20

\vfill\eject


\section{Introduction}
\label{sect-Introduction}
\setcounter{equation}{0} 

The theory of Mean Field Games (MFGs) was introduced independently by Huang-Caines-Malham\'e \cite{HCM06} and Lasry-Lions \cite{LL07a}. Since then, its literature has witnessed a vast increase in various directions and the theory turns out to be extremely rich in applications, including economics \cite{ALLM,LLG}, engineering \cite{HCM071,HCM072}, finance \cite{LackerSoret, LackerZariphopoulou}, social science \cite{BTB,Gelfand} and many others. We refer to Lions \cite{Lions}, Cardaliaguet \cite{Cardaliaguet} and Bensoussan-Frehse-Yam \cite{BFY} for the introduction of the subject in the early stage and Camona-Delarue \cite{CD1, CD2} and Cardaliaguet-Porretta \cite{CP} for more recent developments. Such problems consider limit behavior of large systems where the agents interact with each other in some symmetric way, with the systemic risk as a notable application. The master equation, introduced by Lions \cite{Lions}, characterizes the value of the MFG provided there is a unique mean field equilibrium. It  plays the role of the PDE in the standard literature of controls/games, and is a powerful tool in the mean field framework. The main feature of the master equation is that its state variables include a probability measure $\mu$, representing the distribution of the population, so it can be viewed as a PDE on the Wasserstein space of probability measures. 

In a standard MFG, the interaction is only through the law of the state. In many applications, however, the interaction could be through the joint law of the state and the control. Such a game is called a Mean Field Game of Controls (MFGCs), which was also termed as extended MFGs in the early stage. To be precise, let $B$ and $B^0$ stand for the idiosyncratic and common noises respectively. Given an $\mathbb F^{B^0}$-adapted stochastic measure flow $\{\nu_\cd\}=\{\nu_t\}_{t\in [0,T]}\subset \mathcal{P}_2(\mathbb R^{2d})$, we denote its first marginal by $\mu_t:={\pi_1}_\#\nu_t\in\mathcal{P}_2(\mathbb R^d)$ where $\pi_1(x,a)=x$, for any $(x,a)\in\mathbb R^{d}\times\mathbb R^d$, is a projection. Given above $\{\nu_\cd\}$, we would like to minimize the following cost functional over all admissible controls $\alpha:[0,T]\times\mathbb R^d\times C([0,T];\mathbb R^d) \to\mathbb R^d$: for any $\xi\in \dbL^2(\cF_0, \dbR^d)$, 
\begin{equation}\label{eq:MFGCint}
J(\xi,\{\nu_\cd\};\alpha):=\mathbb E\Big[ G(X_T^{\xi, \{\nu_\cd\},\alpha},\mu_T)+\int_0^Tf(X_t^{\xi,\{\nu_\cd\},\alpha},\alpha(t,X_t^{\xi,\{\nu_\cd\},\alpha},B_{[0,t]}^0),\nu_t)dt\Big],
\end{equation}
where, for a constant $\b\ge 0$,
\begin{equation}\label{eq:Xint}
X_t^{\xi,\{\nu_\cd\},\alpha}=\xi+\int_{0}^tb(X_s^{\xi,\{\nu_\cd\},\alpha},\alpha(s,X_s^{\xi,\{\nu_\cd\},\alpha}, B_{[0,s]}^0),\nu_s)ds+ B_t+\beta B_t^0.
\end{equation}
Here the running drift and cost $b, f$ depend on  the joint law of the state and control, while the terminal cost $G$ depends on the law of the state only. 
We call $(\a^*,\{\nu^*_\cd\})$ a Nash equilibrium if 
\beaa
\alpha^*\in\arg\min_{\alpha}J(\xi,\{\nu_\cd^*\};\alpha),\q\mbox{and}\q \nu_t^*=\mathcal{L}_{(X_t^{\xi, \{\nu^*_\cd\},\alpha^*},\a^*_t)|\cF^{B^0}_t}.
\eeaa
Introduce the Hamiltonian $H$:
\begin{equation}\label{eq:Hint}
H(x,p,\nu):=\inf_{a\in\mathbb R^d}\big[p\cdot b(x,a, \nu)+f(x,a,\nu)\big],~\mbox{with an optimal argument}~ a^* = \phi(x, p, \nu).
\end{equation}
The above problem leads to the following MFGC system of forward-backward stochastic partial differential equations (FBSPDEs) with a solution $(\{\mu_\cd\}, \{\nu_\cd\}, u, v)$:
\begin{equation} 
\label{SPDEint}
\left.\begin{array}{lll}
d\mu_t(x) = \big[\frac{\hat \beta^2}{2} \tr(\partial_{xx} \mu_t(x)) - div\big(\mu_t(x) \partial_p H(x,\partial_x u(t,x),\nu_t)\big)\big]dt-\beta\partial_x\mu_t(x)\cdot d B_t^0;\ms\\
d u(t, x)=  v(t,x)\cdot dB_t^0 - \big[\frac{\hat\beta^2}{2} \tr(\partial_{xx} u(t,x))+\beta\tr(\partial_x v^\top(t,x)) + H(x,\partial_x u(t,x),\nu_t) \big]dt;\ms\\
\nu_t=(id, \phi(\cdot,\pa_x u(t,\cdot),\nu_t))_{\#}\mu_t;\qq \hat \b^2= 1+\b^2; \ms\\ 
\mu_0 =  \cL_\xi,\quad u(T,x) = G(x, \mu_T).
\end{array}\right.
\end{equation}
The wellposedness of the above MFGC system has been investigated by many authors in recent years, essentially in the case $\b=0$ and $b(x,a,\nu) = a$.  For example, Gomes-Patrizi-Voskanyan \cite{GPV}, Kobeissi \cite{K1}, Graber-Mayorga \cite{GraberMayorga} investigated the system under some smallness conditions, and the global wellposedness (especially the uniqueness) was studied by Gomes-Voskanyan \cite{GV0,GV1}, Carmona-Lacker \cite{CarmonaLacker}, Carmona-Delarue \cite{CD1}, Cardaliaguet-Lehalle \cite{CardaliaguetLehalle}, Kobeissi \cite{K0}, under the crucial Lasry-Lions monotonicity condition. We also refer to  Djete \cite{Djete} for some convergence analysis from $N$-player games to  MFGCs and Achdou-Kobeissi \cite{AK} for some numerical studies of MFGCs, without requiring the uniqueness of the equilibria. However, to our best knowledge, the wellposedness of master equations for MFGCs remains completely open. We recall that the master equation is the PDE to characterize the value function $V$ of the MFGC, provided the equilibrium is unique, and it also serves as the decoupling function $V$ of the MFGC system \reff{SPDEint}:
\beaa
u(t,x) = V(t,x, \mu_t).
\eeaa

The monotonicity condition is used to guarantee the uniqueness of the mean field equilibria, and then the global wellposedness of MFG master equations. There are three types of monotonicity conditions in the literature for master equations of standard MFGs: the Lasry-Lions monotonicity, the displacement monotonicity and the  anti-monotonicity.  The Lasry-Lions monotonicity, introduced by Lions \cite{Lions} and extensively used in the literature, can be formulated as following: for any $\xi,\eta\in \dbL^2(\cF_T^1; \dbR^d)$ and their independent copies $\tilde \xi, \tilde \eta$ in the probability space $(\tilde \Omega, \tilde \dbF,\tilde \dbP)$ (see their definitions in Section 2),
\begin{equation}\label{eq:LLint}
\tilde{\mathbb E}\Big[\langle\pa_{x\mu}G(\xi,\mathcal{L}_{\xi},\tilde\xi)\tilde{\eta},\eta\rangle\Big]\geq 0.
\end{equation}
The displacement monotonicity, originating in Ahuja \cite{Ahuja}, is
\begin{equation}\label{eq:disint}
\tilde{\mathbb E}\Big[\langle\pa_{x\mu}G(\xi,\mathcal{L}_{\xi},\tilde\xi)\tilde{\eta},\eta\rangle + \langle\pa_{xx}G(\xi,\mathcal{L}_{\xi})\eta,\eta\rangle\Big]\geq 0,
\end{equation}
which can be further weaken to the displacement semi-monotonicity: for some constant $\l\ge 0$,
\begin{equation}\label{eq:dissemiint}
\tilde{\mathbb E}\Big[\langle\pa_{x\mu}G(\xi,\mathcal{L}_{\xi},\tilde\xi)\tilde{\eta},\eta\rangle + \langle\pa_{xx}G(\xi,\mathcal{L}_{\xi})\eta,\eta\rangle + \l |\eta|^2\Big]\geq 0.
\end{equation}
See, e.g., Bensoussan-Graber-Yam \cite{BGY2} and Gangbo-Meszaros-Mou-Zhang \cite{GMMZ}. It is worth noting that if $G$ is Lasry-Lions monotone and $\pa_{xx}G$ is bounded, then $G$ is displacement semi-monotone. The anti-monotonicity, recently introduced by the authors \cite{MZ3}, takes the following form 
\begin{equation} 
\label{eq:antiint}
\left.\begin{array}{lll}
\qq\,\,\,\,\tilde{\mathbb E}\Big[\l_0\langle\pa_{xx}G(\xi,\mathcal{L}_{\xi})\eta,\eta\rangle+\l_1\langle\pa_{x\mu}G(\xi,\mathcal{L}_{\xi},\tilde \xi)\tilde\eta,\eta\rangle\\
+|\pa_{xx}G(\xi,\mathcal{L}_{\xi})\eta|^2+\l_2\Big|\tilde{\mathbb E}[\pa_{x\mu}G(\xi,\mathcal{L}_{\xi},\tilde\xi)\tilde\eta]\Big|^2-\l_3|\eta|^2\Big]\leq 0,
\end{array}\right.
\end{equation}
for some appropriate constants $\l_0>0$, $\l_1\in\mathbb R$, $\l_2>0$, $\l_3\geq 0$. 

In  \cite{GMMZ, MZ3} we made a simple but crucial observation: the propagation of a monotonicity is crucial for the global wellposedness of the (standard) MFG master equations. That is, provided the terminal condition $G$ satisfies one of the above three types of monotonicity conditions, if one can show a priori that any classical solution $V$ of the master equation satisfies the same type of monotonicity for all time $t$, then  one can establish the global wellposedness of the master equation, which in turn will imply the uniqueness of mean field equilibria and the convergence from the $N$-player game to the MFG. Our goal is to extend all these results to MFGCs, but in this paper we focus only on the propagation of  these three types of monotonicities. That is, we shall follow the approach in \cite{GMMZ, MZ3} to find sufficient conditions on the Hamiltonian $H$ (or alternatively on $b$ and $f$) so that the monotonicity of $G$ can be propagated along $V(t, \cd,\cd)$, provided the master equation has a classical solution $V$. We shall leave the global wellposedness of the master equations and the convergence of the $N$-player games to an accompanying paper.    
 
The Lasry-Lions monotonicity condition has already been used to study the MFGC system \reff{SPDEint}, as mentioned earlier. It is observed in \cite{GMMZ} that, for standard MFGs with non-separable $f$, the Lasry-Lions monotonicity can hardly be propagated. The extra dependence on the law of the control actually helps for propagating the Lasry-Lions monotonicity, in particular, the separability of $f$ is not required anymore. 
 
The displacement semi-monotonicity condition has been introduced in \cite{GMMZ}, however, only the propagation of displacement monotonicity is established there. In this paper, we manage to propagate the displacement semi-monotonicity for MFGCs, so it improves the result of \cite{GMMZ} even for standard MFGs.  In particular, by combining with the arguments in \cite{GMMZ}, we obtain easily the global wellposedness result of standard MFG master equations under displacement semi-monotonicity conditions. We remark again that the displacement semi-monotonicity is weaker than both displacement monotonicity and Lasry-Lions monotonicity (provided $\pa_{xx} V$ is bounded, which is typically the case), so in this sense our result  provides a unified framework for the wellposedness theory of master equations under the Lasry-Lions monotonicity and displacement monotonicity conditions.
 
Another feature of our results is that we allow for a general form of the drift $b$. In the literature, one typically sets $b(\cd, a, \cd) = a$ (or slightly more general forms), and then focuses on appropriate monotonicity conditions of $f$ to ensure the uniqueness of the mean field equilibria and/or  the wellposedness of the master equations. However, for a general $b$, especially when $b$ depends on the law (of the state and/or the control), it does not make sense to propose monotonicity condition on $f$ alone.  A conceivable notion of monotonicity on the general $b$ has never been studied, to our best knowledge. Our approach works on the Hamiltonian $H$ directly, which has the mixed impacts of $b$ and $f$ together. Again, our results are new in this aspect even for standard MFGs.

The rest of the paper is organized as follows. In Section \ref{sect-MFGC} we introduce MFGCs.  In Section \ref{sect-master} we introduce the master equation and the notions of monotonicities. In Sections \ref{sect-LL}, \ref{sect-displacement}, and \ref{sect-anti} we propagate the three types of monotonicities, one in each section.  In particular, in Subsection \ref{sect-wellposedness} we also establish the global wellposedness of standard MFG master equations under displacement semi-monotonicity conditions. Finally, some technical proofs are postponed to Appendix.

\section{Mean Field Games of Controls}
\label{sect-MFGC} 
\setcounter{equation}{0}
We consider the setting in \cite{GMMZ}. Let $d$ be a dimension and $[0, T]$ a fixed finite time horizon. Let $(\O_0, \dbF^0, \dbP_0)$ and $(\O_1, \dbF^1, \dbP_1)$ be two filtered probability spaces, on which are defined $d$-dimensional Brownian motions $B^0$ and $B$, respectively. For $\dbF^i =\{\cF^i_t\}_{0\le t\le T}$, $i=0,1$, we assume $\cF^0_t=\cF^{B^0}_t$,  $\cF^1_t =\cF^1_0 \vee \cF^{B}_t$, and $\dbP_1$ has no atom in $\cF^1_0$ so it can support any measure on $\dbR^d$ with finite second order moment. Consider the product spaces 
\begin{equation}\label{product}
\O := \O_0 \times \O_1,\q \dbF = \{\cF_t\}_{0\le t\le T} := \{\cF^0_t \otimes \cF^1_t\}_{0\le t\le T},\q \dbP := \dbP_0\otimes \dbP_1,\q \dbE:= \dbE^\dbP.
\end{equation}
In particular, $\cF_t := \si(A_0\times A_1: A_0\in \cF^0_t, A_1\in \cF^1_t\}$ and $\dbP(A_0\times A_1) = \dbP_0(A_0) \dbP_1(A_1)$. 
 We shall automatically extend $B^0, B, \dbF^0, \dbF^1$ to the product space in the obvious sense, but using the same notation. Note that $B^0$ and $B^1$ are independent $\dbP$-Brownian motions and are independent of $\cF_0$. 

It is convenient to introduce another filtered probability space $(\tilde \O_1, \tilde \dbF^1, \tilde B, \tilde \dbP_1)$  in the same manner as $(\O_1,  \dbF^1, B, \dbP_1)$, and consider the larger filtered probability space given by
\begin{equation}\label{product2}
\tilde \O := \O\times \tilde\O_1 ,\q \tilde\dbF = \{\tilde \cF_t\}_{0\le t\le T} := \{\cF_t \otimes \tilde \cF^1_t\}_{0\le t\le T},\q \tilde \dbP := \dbP\otimes \tilde\dbP_1,\q \tilde \dbE:= \dbE^{\tilde \dbP}.
\end{equation}
Given an $\cF_t$-measurable random variable $\xi = \xi(\o^0, \o^1)$, we say $\tilde \xi = \tilde \xi(\o^0, \tilde \o^1)$ is a conditionally independent copy of $\xi$ if, for each $\o^0$, the $\dbP_1$-distribution of $\xi(\o^0, \cd)$ is equal to the $\tilde\dbP_1$-distribution of $\tilde\xi(\o^0, \cd)$.  That is, conditional on $\cF^0_t$, by extending to $\tilde \O$ the random variables $\xi$ and $\tilde \xi$ are conditionally independent and have the same conditional distribution under $\tilde \dbP$. Note that, for any appropriate deterministic function $\f$,
\bea
\label{conditional expectation}
\left.\ba{c}
\dis \tilde \dbE_{\cF^0_t} \big[ \f(\xi, \tilde \xi)\big] (\o^0) = \dbE^{\dbP_1\otimes \tilde\dbP_1}\Big[\f\big(\xi(\o^0, \cdot), \tilde\xi(\o^0, \tilde\cdot)\big)\Big],\q \dbP_0-\mbox{a.e.}~\o^0;\\
\dis \tilde \dbE_{\cF_t} \big[ \f(\xi, \tilde \xi)\big] (\o^0,\o^1) = \dbE^{\tilde\dbP_1}\Big[\f\big(\xi(\o^0, \o^1), \tilde\xi(\o^0, \tilde \cdot)\big)\Big],\q \dbP-\mbox{a.e.}~(\o^0, \o^1).
\ea\right.
\eea
Here $ \dbE^{\tilde\dbP_1}$ is the expectation on $\tilde \o^1$, and $\dbE^{\dbP_1\times \tilde\dbP_1}$ is on $(\o^1, \tilde \o^1)$. Throughout the paper, we will use the probability space $(\O, \dbF, \dbP)$. However, when conditionally independent copies of random variables or processes are needed, we will tacitly use the extension to the larger space $(\tilde \O, \tilde \dbF, \tilde \dbP)$ without mentioning.

When we need two conditionally independent copies, we introduce further $(\bar \O_1, \bar \dbF^1, \bar B, \bar \dbP_1)$ and the product space $(\bar\O, \bar \dbF, \bar\dbP, \bar\dbE)$  as in \reff{product2}, and set the joint product space
\begin{align}\label{product3}
\bar{\tilde \O}  := \O\times \tilde\O_1\times \bar\O_1 ,~ \bar{\tilde \dbF} = \{\bar{\tilde \cF}_t\}_{0\le t\le T} := \{\cF_t\otimes\tilde\cF^1_t \otimes \bar \cF^1_t\}_{0\le t\le T}, ~\bar {\tilde \dbP} := \dbP\otimes \tilde\dbP_1\otimes \bar\dbP_1,~ \bar{\tilde \dbE}:= \dbE^{\bar{\tilde \dbP}}.
\end{align}

For any dimension $k$ and any constant $p\ge 1$, let $\cP(\dbR^k)$ denote the set of probability measures on $\dbR^k$, and $\cP_p(\dbR^k)$ the subset of $\mu\in \cP(\dbR^k)$ with finite $p$-th moment, equipped with the $p$-Wasserstein distance $W_p$. Moreover, for any sub-$\sigma$-algebra $\mathcal{G}\subset \mathcal{F}_T$, $\mathbb L^p(\mathcal{G})$ denotes the set of $\mathbb R^k$-valued, $\mathcal{G}$-measurable, and $p$-integrable random variables; and for any $\mu\in \cP_p(\dbR^k)$, $\mathbb L^p(\mathcal{G};\mu)$ denotes the set of $\xi\in\mathbb L^p(\mathcal{G})$ with law $\mathcal{L}_{\xi}=\mu$. Similarly, for any sub-filtration $\dbG \subset \dbF$, $\dbL(\dbG; \dbR^k)$ denotes the set of $\dbG$-progressively measurable $\dbR^k$-valued processes.

For a continuous function $U:\cP_2(\dbR^k)\to\dbR$, we recall its linear functional derivative ${\d U\over \d\mu}: \cP_2(\dbR^k)\times\dbR^k\to\dbR$ and Lions derivative $\pa_\mu U: \cP_2(\dbR^k)\times\dbR^k\to\dbR^k$. 
We say $U\in\cC^1(\cP_2(\dbR^k))$ if  $\pa_\mu U$ exists and is continuous on $\cP_2(\dbR^k)\times\dbR^k$, and we note that $\pa_\mu U(\mu, \tilde x) = \pa_{\tilde x} {\d U\over \d\mu}(\mu, \tilde x)$.  Similarly we can define the second order derivative $\pa_{\mu\mu} U(\mu,\tilde x,\bar x)$, and we say $U\in \cC^2(\cP_2(\dbR^k))$ if $\pa_\mu U$, $\pa_{\tilde x\mu}U$ and $\pa_{\mu\mu}U$ exist and are continuous. We refer to \cite[Chapter 5]{CD1} or \cite{GT} for more details.

Our mean field game of controls (MFGC) will depend on the following data: 
\beaa
b: \mathbb R^{2d}\times\mathcal{P}_{2}(\mathbb R^{2d})\to\mathbb \dbR^d;\q f: \mathbb R^{2d}\times\mathcal{P}_2(\mathbb R^{2d})\to\mathbb \dbR;\q G: \mathbb R^{d}\times\mathcal{P}_2(\mathbb R^{d})\to\mathbb R;~\mbox{and}~ \b \in [0, \infty).
\eeaa
We shall always assume appropriate technical conditions so that all the equations in this section are wellposed and all the involved random variables are integrable. Given $t_0\in [0,T]$, denote $B_t^{t_0}:=B_t-B_{t_0}$, $B_t^{0,t_0}:=B_t^0-B_{t_0}^0$, $t\in [t_0,T]$.  Let $\cA_{t_0}$ denote the set of admissible controls $\a: [t_0,T]\times \mathbb R^{d}\times C([t_0,T];\mathbb R^d)\to\mathbb R^d$ which are progressively measurable and adapted in the path variable and square integrable; and $\dbL^2(\dbF^{B^{0, t_0}}; \cP_2(\dbR^{2d}))$ the set of $\dbF^{B^{0, t_0}}$-progressively measurable stochastic measure flows $\{\nu_\cdot\}=\{\nu_t\}_{t\in [t_0,T]}\subset \cP_2(\mathbb R^{2d})$.  Here for notational simplicity we assume the controls also take values in $\dbR^d$, and $b$ and $f$ do not depend on time, but one can remove these constraints without any difficulty.

Given $t_0\in [0,T]$, $x\in\mathbb R^d$, $\alpha\in\mathcal{A}_{t_0}$, and $\{\nu_\cd\}\in \dbL^2(\dbF^{B^{0, t_0}}; \cP_2(\dbR^{2d}))$, 
the state of the agent satisfies the following controlled SDE on $[t_0,T]$: 
\bea
\label{Xx}
\left.\ba{c}
\dis X_t^{ \{\nu_\cdot\},\a}=x+\int_{t_0}^{t}b(X_s^{ \{\nu_\cdot\},\a},\alpha_s, \nu_s)ds+B_t^{t_0}+\beta B_t^{0,t_0};\\
\dis \mbox{where}\q X^{ \{\nu_\cdot\},\a} = X^{t_0, \{\nu_{\cdot}\}; x, \alpha},\q \a_s := \alpha (s,X_s^{ \{\nu_\cdot\},\a}, B^{0,t_0}_{[t_0,s]}).
\ea\right.
\eea
Consider the  expected cost for the MFGC: denoting by ${\pi_1}_{\#} \nu_{t} $ the first component of $\nu_{t}$,
\bea
\label{Ja}
J(t_0,x; \{\nu_\cdot\},\a):= \inf_{\alpha\in\mathcal{A}_{t_0}}\mathbb E\Big[G(X_{T}^{ \{\nu_\cdot\},\a},{\pi_1}_{\#} \nu_T)
+\int_{t_0}^Tf(X_t^{ \{\nu_\cdot\},\a}, \alpha_t, \nu_t)dt\Big].
\eea

\begin{defn}
For any $(t,\mu)\in [0,T]\times\mathcal{P}_2(\mathbb R^d)$,   we say $(\alpha^*,\{\nu_\cdot^*\})\in \cA_t \times \dbL^2(\dbF^{B^{0, t_0}}; \cP_2(\dbR^{2d}))$ is a mean field equilibrium (MFE) at $(t,\mu)$  if 
\bea\label{MFE}
\left.\ba{c}
\dis J(t,x;\{\nu_\cd^*\},\alpha^*) = \inf_{\alpha\in\mathcal{A}_{t}}J(t,x;\{\nu_\cd^*\},\alpha),\q \text{for $\mu$-a.e. $x\in\mathbb R^d$};\\
\dis \pi_{1\#}\nu_t^* =\mu,\q \nu_s^*:=\mathcal{L}_{(X_s^*,\alpha^* (s,X_s^*, B^{0,t_0}_{[t_0,s]}))|\cF^0_s}, \q\mbox{where}\\
\dis X_t^*=\xi+\int_{t_0}^tb(X_s^*,\alpha^* (s,X_s^{*}, B^{0,t_0}_{[t_0,s]}),\nu_s^*)ds+B_t^{t_0}+\beta B_t^{0,t_0},\q \xi\in \dbL^2(\cF^1_t, \mu).
\ea\right.
\eea
\end{defn}

When there is a unique MFE for each $(t,\mu)\in [0,T]\times\mathcal{P}_2(\mathbb R^d)$, denoted as $(\alpha^*(t,\mu; \cd),\{\nu_\cd^*(t,\mu)\})$, then the game problem leads to the following value function for the agent: 
\begin{equation}\label{value}
V(t,x,\mu):=J(t,x;\{\nu_\cd^*(t,\mu)\}, \alpha^*(t,\mu; \cd))\quad\text{for any $x\in\mathbb R^d$}.
\end{equation}
We remark that, by \reff{MFE} the above $V$ is well defined only for $\mu$-a.e. $x$. However, for each $t$, its continuous extension to $\dbR^d\times \cP_2(\dbR^d)$ is unique,  and we shall always consider this continuous extension. 
Our goal is to study the master equation for the value function $V(t,x,\mu)$. 

For this purpose, we introduce the Hamiltonian: for $(x, p, \nu)\in \dbR^d\times \dbR^d\times \cP_2(\dbR^{2d})$, 
\bea
\label{H}
H(x,p,\nu):= \inf_{a\in \dbR^d} h(x,p, \nu, a),\q  h(x,p, \nu, a):= p\cdot b(x,a,\nu)+f(x,a,\nu).
\eea
Note that $H$ depends on $\nu$, while $V$ depends only on $\mu = {\pi_1}_{\#} \nu$. We also remark that the Hamiltonian in \cite{GMMZ, MZ3} is $-H$.  To introduce the master equation, which we will do in the next section, we need the following fixed point.

\begin{assum}
\label{assum-fix}
(i) The Hamiltonian $H$ has a unique minimizer $a^* = \phi(x, p, \nu)$, namely
\bea
\label{I}
H(x,p,\nu) = h(x,p, \nu, \phi(x, p, \nu)).
\eea

(ii) For any $\xi \in \dbL^2(\cF)$ and $\eta \in \dbL^2(\si(\xi))$, the following mapping on $\cP_2(\dbR^{2d})$: 
\bea
\label{cI}
\cI^{\xi,\eta}(\nu) := \cL_{(\xi, \phi(\xi, \eta, \nu))}
\eea
has a unique fixed point $\nu^*$: $\cI^{\xi,\eta}(\nu^*)=\nu^*$, and we shall denote  it as  $\Phi(\cL_{(\xi, \eta)})$. 
\end{assum}
We refer to \cite[Lemma 4.60]{CD1} for some sufficient conditions on the existence of $\Phi$. 
By \reff{I} one can  easily check that
\bea
\label{Hp}
b(x,\phi(x,p,\nu),\nu)=\pa_pH(x,p,\nu), \q  f(x,\phi(x,p,\nu),\nu)= H(x, p, \nu) - p\cd \pa_pH(x,p,\nu).
\eea
As in the standard MFG theory, provided $V$ is smooth, $p$ corresponds to $\pa_x V(t, x, \mu)$. Consequently,  later on the above fixed point will be applied as follows: given $(t, \mu)$ and $\xi \in \dbL^2(\cF^1_t, \mu)$, 
\bea
\label{optimal}
\eta = \pa_x V(t, \xi, \mu),\q \nu^* := \Phi(\cL_{\xi, \pa_x V(t, \xi, \mu)}),\q \a^* := \phi(\xi, \pa_x V(t, \xi, \mu), \Phi(\cL_{\xi, \pa_x V(t, \xi, \mu)})).
\eea
Pluging these into  \reff{MFE} we obtain the following McKean-Vlasov SDE: recalling \reff{Hp},
\bea
\label{X*}
\left.\ba{c}
\dis X_t^{*}=\xi+\int_{t_0}^{t} \pa_p H\big(X_s^{*}, \pa_x V(s, X^*_s, \mu^*_s\big), \nu^*_s)ds+B_t^{t_0}+\beta B_t^{0,t_0},\\
\dis\mbox{where}\q  \mu^*_s := \cL_{X^*_s|\cF^0_s},\q \nu^*_s := \Phi\big(\cL_{(X^*_s, \pa_x V(s, X^*_s, \mu^*_s))|\cF^0_s}\big). 
\ea\right.
\eea
That is, if $V$ is smooth, then under Assumption \ref{assum-fix} we may obtain the unique MFE $\a^*$ through \reff{optimal} and \reff{X*} (by abusing the notation $\a^*$): given $(t_0, \mu)$ and $\xi\in \dbL^2(\cF^1_{t_0}, \mu)$,
\bea
\label{a*}
\a^*(s, x, B^{0, t_0}_{[t_0, s]}) =  \phi(x, \pa_x V(s, x, \mu^*_s), \nu^*_s).
\eea
Here we used the fact that $\mu^*_s, \nu^*_s$ are actually adapted to the shifted filtration generated by $B^{0, t_0}$.

Assumption \ref{assum-fix} (i) is more or less standard in the literature, for example when $h$ in \reff{H} is convex in $a$. In particular, when $b(x,a,\nu) = a$, which is often the case in the literature, we have $\phi = \pa_p H$. We next provide two examples for Assumption \ref{assum-fix} (ii). 
\begin{eg}
\label{eg-separable}
Assume $b, f$ are separable in the following sense:
\bea
\label{separable}
b(x, a, \nu) = b_0(x, a, {\pi_1}_{\#} \nu) + b_1(x, \nu),\q  f(x, a, \nu) = f_0(x, a, {\pi_1}_{\#} \nu) + f_1(x, \nu).
\eea
In this case \reff{H} becomes:
\beaa
&\dis H(x,p,\nu) = H_0(x, p,  {\pi_1}_{\#} \nu) + H_1(x,  p, \nu),\q \mbox{where}\\
&\dis H_0(x,p,\mu):=\inf_{a\in \dbR^d} \big[p\cdot b_0(x,a,\mu)+f_0(x,a,\mu)\big],\q H_1(x,p, \nu):= p\cd b_1(x, \nu) + f_1(x, \nu).
\eeaa
Assume Assumption \ref{assum-fix} (i) holds, and clearly in this case we have $a^* = \phi(x, p, \mu)$, with the dependence on $\nu$ only through its first component $\mu={\pi_1}_{\#} \nu$.  Then $\cI^{\xi,\eta}(\nu) := \cL_{(\xi, \phi(\xi, \eta, {\pi_1}_{\#} \nu))}$.  Notice further that the fixed point requires ${\pi_1}_{\#} \nu^* = \cL_\xi$. Then it is obvious that Assumption \ref{assum-fix} (ii) holds with $\Phi(\cL_{(\xi, \eta)}) = \cL_{(\xi, \phi(\xi, \eta, \cL_\xi))}$.
\end{eg}

We note that the above $f$ satisfies the conditions in \cite[Lemma 4.60]{CD1}, while the drift $b$ is more general. The next example, however,  is out of the scope of \cite[Lemma 4.60]{CD1}.
\begin{eg}
\label{eg-nonseparable}
Assume $d=1$ and, by writing $\dbE_\nu[\a]$ to indicate expectation under law $\cL_{\a} = {\pi_2}_{\#} \nu$,
\bea
\label{nonseparable}
b(x, a, \nu) = -b_0(x, {\pi_1}_{\#} \nu) a + b_1(x, \nu), \q f(x, a, \nu) = {|a|^2\over 2} - a f_0(x, {\pi_1}_{\#} \nu, \dbE_{\nu}[\a]) + f_1(x, \nu).
\eea
One can easily see that $\phi(x, p, \nu) = f_0(x, {\pi_1}_{\#} \nu, \dbE_\nu[\a])+p b_0(x, {\pi_1}_{\#} \nu)$ and thus
\beaa
\cI^{\xi,\eta}(\nu) := \cL_{(\xi, ~ f_0(\xi, {\pi_1}_{\#} \nu,  \dbE_\nu[\a])+ b_0(\xi, {\pi_1}_{\#} \nu)\eta)}.
\eeaa
Then $\cI^{\xi,\eta}$ has a fixed point if and only if the following mapping has a fixed point:
\bea
\label{psi}
m\in \dbR \to \psi^{\xi, \eta}(m):= \dbE\big[ f_0(\xi, \cL_\xi,  m)+ b_0(\xi, \cL_\xi)\eta\Big].
\eea
Assume $\pa_m f_0 \le 1-\e$ for some $\e>0$, in particular if $f_0$ is decreasing in $m$, then $\pa_m \psi^{\xi, \eta}\le 1-\e$ and thus $\psi^{\xi, \eta}$ has a unique fixed point $m^* = \f(\cL_{(\xi, \eta)})$. Therefore, $\cI^{\xi,\eta}$ has a unique fixed point:
\bea
\label{fixpsi}
\Phi(\cL_{(\xi, \eta)}) =  \cL_{(\xi, ~ f_0(\xi, \cL_\xi,  \f(\cL_{(\xi, \eta)}))+ b_0(\xi, \cL_\xi)\eta)}.
\eea
\end{eg}

\subsection{Derivatives of measure valued functions}
Note that $\Phi$ is a mapping from $\cP_2(\dbR^{2d})$ to $\cP_2(\dbR^{2d})$. Consider an arbitrary dimension $k$. In this subsection we introduce the linear functional derivative of functions mapping from $\cP_2(\dbR^k)$ to $\cP_2(\dbR^k)$, which is interesting in its own right. We refer to \cite[Eq. (5.52)]{CD1} for the linear functional derivative of functions mapping from $\cP_2(\dbR^k)$ to $\dbR$.  Let $\mathcal{S}(\mathbb R^k)$ denote the Schwartz space, namely the set of smooth functions $u\in \cC^\infty(\dbR^k; \dbR)$ such that $u$ and all its derivatives decrease rapidly when $|x|\to \infty$, and $\cS'(\dbR^k)$ its dual space, namely the space of tempered distributions.

\begin{defn}
\label{defn-derivative}
Consider a mapping $\Phi: \cP_2(\dbR^k) \to \cP_2(\dbR^k)$. We say ${\d \Phi\over \d\rho}: \cP_2(\dbR^k) \times \dbR^k \to \mathcal{S}'(\mathbb R^k)$ is the linear functional derivative of $\Phi$ if,  for any $\psi \in \mathcal{S}(\mathbb R^k)$, 
\bea
\label{derivative}
\frac{\delta \Psi}{\delta\rho} (\rho, x) = \Big\la {\d \Phi\over \d\rho}(\rho, x), ~ \psi\Big\ra,\q\mbox{where}\q \Psi(\rho) :=  \int_{\dbR^k} \psi(x) \Phi(\rho; dx).
\eea
\end{defn}
We note that $ {\d \Phi\over \d\rho}(\rho, x)$ is well defined for $\rho$-a.e. $x$.  

For the applications later, we will require ${\d \Phi\over \d\rho}$ to have stronger properties. For this purpose,  let $\mathcal{SM}_2(\mathbb R^k)$ denote the set of the square integrable signed measures of bounded variation on $\mathbb R^k$. That is, $m$ has the unique decomposition $m=m_1-m_2$ and $\int_{\mathbb R^k}(1+|y|^2)|m|(dy)<\infty$, where $m_1,m_2$ are mutually singular non-negative measures on $\mathbb R^k$, and $|m|(dy):=m_1(dy)+m_2(dy)$, see e.g. \cite{Bogachev} for details. Moreover, for any $n\ge 0$, let $\mathcal{DSM}^n_2(\dbR^k)\subset \cS'(\dbR^k)$ denote the linear span of generalized derivatives of signed measures in $\mathcal{SM}_2(\mathbb R^k)$ up to order $n$, namely the span of terms taking the form $\pa_{y_1}^{j_1}\cds \pa_{y_k}^{j_k} m$, where $m\in \mathcal{SM}_2(\mathbb R^k)$ and $\sum_{i=1}^k j_i \le n$. On the other hand, let $\cC^n_2(\dbR^k)$ denote the set of functions $\psi: \dbR^k \to \dbR$ such that $\psi$ has  continuous derivatives up to order $n$ and 
\bea
\label{n-norm}
\|\psi\|_n := \sup_{y\in \dbR^k} \sum_{j_1+\cds+j_k\le n} {|\pa_{y_1}^{j_1}\cds \pa_{y_k}^{j_k} \psi(y)|\over 1+|y|^2} <\infty.
\eea
Then clearly $\mathcal{DSM}^n_2(\dbR^k)$ is in the dual space of $\cC^n_2(\dbR^k)$ in the sense that
\bea
\label{integration-by-parts}
\Big\la \pa_{x_1}^{j_1}\cds \pa_{x_k}^{j_k} m,~ \psi \Big\ra = (-1)^{\sum_{i=1}^k j_i} \int_{\dbR^k} \pa_{y_1}^{j_1}\cds \pa_{y_k}^{j_k} \psi(y) m(dy).
\eea
Now if ${\d \Phi\over \d\rho}(\rho, x) \in \mathcal{DSM}^n_2(\dbR^k)$, then we may extend \reff{derivative} to all $\psi \in \cC^n_2(\dbR^k)$, and we shall write
\beaa
\int_{\dbR^k} \psi(y) {\d \Phi\over \d\rho}(\rho, x; dy):= \big\la {\d \Phi\over \d\rho}(\rho, x), ~\psi\big\ra,\q \forall \psi\in \cC^n_2(\dbR^k),
\eeaa
where the right side is in the sense of \reff{integration-by-parts}.

We next show two examples. 
\begin{eg}
\label{eg-trivial}
Let $\Phi(\rho)=\rho$ for any $\rho\in\cP_2(\mathbb R^k)$. Then $\frac{\delta\Phi}{\delta\rho}(\rho,x; dy)=\delta_x(dy)$, namely $ \frac{\delta\Phi}{\delta\rho}(\rho,x) \in {\cal SM}_2(\dbR^k) = {\cal DSM}^0_2(\dbR^k)$ for all  $\rho\in\mathcal{P}_2(\mathbb R^d)$ and $x\in\mathbb R^k$.
\end{eg}
\proof For any $\psi\in\mathcal{S}(\mathbb R^k)$, by \eqref{derivative} we have $\Psi(\rho) = \int_{\dbR^d} \psi(x) \rho(dx)$. Then
$
\frac{\delta \Psi}{\delta\rho}(\rho,x)=\psi(x)=\int_{\mathbb R^k}\psi(y)\delta_x(dy),
$
and thus  $\frac{\delta\Phi}{\delta\rho}(\rho,x)=\delta_x\in\mathcal{SM}_2(\dbR^k)$. 
\qed

\begin{eg}
\label{eg-derivative}
Set $\Phi(\cL_{(\xi, \eta)}) := \cL_{(\xi, \eta + c\dbE[\eta])}$, $\forall \xi, \eta\in \dbL^2(\cF; \dbR^d)$, for some constant $c\in \dbR$. Then $ \frac{\delta\Phi}{\delta\rho}(\rho,x, p) \in  {\cal DSM}^1_2(\dbR^{2d})$. More precisely, letting $\dbE_\rho$ denote expectation under law  $\rho=\cL_{(\xi, \eta)}$,
\bea
\label{derivative-example}
\dis {\d \Phi\over \d\rho}(\rho, x, p; d\tilde x, d\tilde p) = \d_x (d\tilde x) \d_{p+ c\dbE_\rho[\eta]}(d\tilde p)    - c ~\pa_{\tilde p}\Phi(\rho)(d\tilde x, d\tilde p)\cdot p.
\eea
\end{eg}
\proof For any $\psi \in \cS(\dbR^{2d})$, we have $\Psi(\rho)= \dbE_{\rho}\big[\psi(\xi, \eta + c\dbE_{\rho}[\eta])\big]$. Then
\beaa
\frac{\delta\Psi}{\delta\rho}(\rho,x,p) &=& \psi\big(x,p+c\dbE_{\rho}[\eta]\big)+ c\dbE_{\rho}\big[\pa_p\psi(\xi, \eta + c\dbE_{\rho}[\eta])\big]\cdot p\\
&=& \psi\big(x,p+c\dbE_{\rho}[\eta]\big)+ cp \cdot \int_{\dbR^k} \pa_{\tilde p} \psi(\tilde x, \tilde p) \Phi(\rho)(d\tilde x, d \tilde p).
\eeaa
Compare this with \reff{derivative}, we obtain \reff{derivative-example} immediately.
\qed

Our main result of this part is the following chain rule. We shall use the notation $\nu = \Phi(\rho)$.
\begin{prop}\label{prop-chain}
Let $\Phi: \cP_2(\dbR^k)\to \cP_2(\dbR^k)$, $U: \cP_2(\dbR^k)\to \dbR$. Assume  

(i)  $\Phi$ has a linear functional derivative $\frac{\delta\Phi}{\delta\rho}(\rho,x)\in\mathcal{DSM}^n_2(\dbR^k)$ for all $(\rho,x)\in\cP_2(\dbR^k)\times\dbR^k$;  $\frac{\delta\Phi}{\delta\rho}(\rho,x)$ is continuous in $(\rho, x)$ under the weak topology, that is, for any $\psi \in \cC^n_2(\dbR^k)$, the mapping $(\rho, x)\to  \big\la\frac{\delta\Phi}{\delta\rho}(\rho,x),  \psi\big\ra$ is continuous (under $\cW_2$ for $\rho$); and, for any compact set $K \subset \cP_2(\dbR^k)$, there exists a constant $C_K>0$ such that
\bea
\label{Phiquadratic}
\sup_{\rho\in K}\Big|\big\la \frac{\delta\Phi}{\delta\rho}(\rho, x),~ \psi\big\ra\Big| \le C_K \|\psi\|_n [1+|x|^2],\q \forall \psi \in \cC^n_2(\dbR^k).
\eea

(ii)  $U$ has a linear functional derivative $\frac{\delta U}{\delta \nu}$;  for each $\nu \in \cP_2(\dbR^k)$, $\frac{\delta U}{\delta\nu}(\nu, \cd) \in \cC^n_2(\dbR^k)$;  and, by equipping $\cC^n_2(\dbR^k)$ with the norm  $\|\cd\|_n$ in \reff{n-norm},  the mapping $\nu \to \frac{\delta U}{\delta\nu}(\nu, \cd)$ is continuous.
 
 Then the composite function $\wh U:= U\circ \Phi:\cP_2(\dbR^k)\to\mathbb R$ has a linear functional derivative:
\bea
\label{linearchain}
\frac{\delta \wh U}{\delta\rho} (\rho, x) =  \int_{\dbR^k} \frac{\delta U}{\delta \nu} (\Phi(\rho), y) {\d \Phi\over \d\rho}(\rho, x; dy).
\eea
\end{prop}
\proof Fix $\rho, \rho' \in \cP_2(\dbR^k)$. For $0< \e<1$, denote $\rho_\e:= \rho + \e(\rho'-\rho)$. By the definition of ${\d U\over \d\nu}$ we have
\beaa
&\dis \wh U(\rho_\e) - \wh U(\rho) = U(\Phi(\rho_\e)) - U(\Phi(\rho)) =\int_0^1 \Big[\Psi_\th(\rho_\e) - \Psi_\th(\rho)\Big] d\th,\\
&\dis \mbox{where}\q \psi_\th(x) := {\d U\over \d \nu}\Big(\th \Phi(\rho_\e) +(1-\th) \Phi(\rho), x\Big),\q \Psi_\th(\tilde\rho) := \int_{\dbR^k} \psi_\th(x) \Phi(\tilde\rho; dx),~ \forall \tilde \rho\in \cP_2(\dbR^k).
\eeaa
 Then,  by \reff{derivative} we have
 \beaa
 \frac{\delta  \Psi_\th}{\delta\rho}(\tilde \rho, x) =  \int_{\dbR^k}  \psi_\th(y) {\d \Phi\over \d\rho}(\tilde\rho, x; dy) =  \int_{\dbR^k} \frac{\delta U}{\delta\nu} \Big(\th \Phi(\rho_\e) +(1-\th) \Phi(\rho), y\Big) {\d \Phi\over \d\rho}(\tilde\rho, x; dy).
 \eeaa
Note that $\rho + \tilde \th (\rho_\e - \rho) = \rho_{\tilde \th \e}$, then
\beaa
\left.\ba{lll}
\dis {1\over \e}\big[ \wh U(\rho_\e) - \wh U(\rho) \big] = {1\over \e} \int_0^1 \int_0^1 \int_{\dbR^k}{\d \Psi_\th\over \d \rho} (\rho_{\tilde \th\e}, x) (\rho_\e - \rho)(dx) d\tilde \th d\th\ms\\
\dis = \int_0^1 \int_0^1 \int_{\dbR^k}  \int_{\dbR^k} \frac{\delta U}{\delta\nu} \Big(\th \Phi(\rho_\e) +(1-\th) \Phi(\rho), y\Big) {\d \Phi\over \d\rho}(\rho_{\tilde \th \e}, x; dy) (\rho' - \rho)(dx) d\tilde \th d\th\ms\\
\dis= I_1(\e) + I_2(\e),
\ea\right.
\eeaa
where
\beaa
I_1(\e) &:=& \int_0^1 \int_{\dbR^k}  \int_{\dbR^k} \frac{\delta U}{\delta\nu} \Big(\Phi(\rho), y\Big) {\d \Phi\over \d\rho}(\rho_{\tilde \th \e}, x; dy) (\rho' - \rho)(dx) d\tilde \th;\\
I_2(\e) &:=&  \int_0^1 \int_0^1 \int_{\dbR^k}  \int_{\dbR^k} \Big[\frac{\delta U}{\delta\nu} \Big(\th \Phi(\rho_\e) +(1-\th) \Phi(\rho), y\Big)- \frac{\delta U}{\delta\nu} \big(\Phi(\rho), y\big)\Big] \times\\
&& {\d \Phi\over \d\rho}(\rho_{\tilde \th \e}, x; dy) (\rho' - \rho)(dx) d\tilde \th d\th.\
\eeaa
Clearly $\lim_{\e\to 0}W_2(\rho_{\tilde \th \e}, \rho)=0$. By the continuity of ${\d \Phi\over \d\rho}$ we have
\beaa
\lim_{\e\to 0}\int_{\dbR^k} \frac{\delta U}{\delta\nu} \Big(\Phi(\rho), y\Big) {\d \Phi\over \d\rho}(\rho_{\tilde \th \e}, x; dy)=\int_{\dbR^k} \frac{\delta U}{\delta\nu} \Big(\Phi(\rho), y\Big) {\d \Phi\over \d\rho}(\rho, x; dy),\q \forall \tilde \th, x.
\eeaa
Moreover, note that $K:= \{\rho_\e: 0\le \e\le 1\} \subset \cP_2(\dbR^k)$ is  compact. Then by \reff{Phiquadratic} we have
\beaa
\Big|\int_{\dbR^k} \frac{\delta U}{\delta\nu} \Big(\Phi(\rho), y\Big) {\d \Phi\over \d\rho}(\rho_{\tilde \th \e}, x; dy)\Big| \le C\|\frac{\delta U}{\delta\nu} \big(\Phi(\rho), \cd\big)\|_n [1+|x|^2].
\eeaa
Now it follows from the dominated convergence theorem that
\bea
\label{I1e}
\lim_{\e\to 0} I_1(\e) =\int_{\dbR^k}\int_{\dbR^k} \frac{\delta U}{\delta\nu} \Big(\Phi(\rho), y\Big) {\d \Phi\over \d\rho}(\rho, x; dy) (\rho'-\rho)(dx).
\eea

Moreover,  by \reff{Phiquadratic} again we have
\beaa
&&\dis \Big|\int_{\dbR^k} \Big[\frac{\delta U}{\delta\nu} \Big(\th \Phi(\rho_\e) +(1-\th) \Phi(\rho), y\Big)- \frac{\delta U}{\delta\nu} \big( \Phi(\rho), y\big)\Big]  {\d \Phi\over \d\rho}(\rho_{\tilde \th \e}, x; dy) \Big|\\
&&\dis\le C\|\frac{\delta U}{\delta\nu} \Big(\th \Phi(\rho_\e) +(1-\th) \Phi(\rho), \cd\Big)- \frac{\delta U}{\delta\nu} \big( \Phi(\rho), \cd\big)\|_n [1+|x|^2].
\eeaa
Then
\beaa
|I_2(\e)| &\le& C\int_0^1 \int_{\dbR^k}\|\frac{\delta U}{\delta\nu} \Big(\th \Phi(\rho_\e) +(1-\th) \Phi(\rho), \cd\Big)- \frac{\delta U}{\delta\nu} \big( \Phi(\rho), \cd\big)\|_n [1+|x|^2] (\rho'+\rho)(dx) d\th\\
&\le& C\int_0^1 \|\frac{\delta U}{\delta\nu} \Big(\th \Phi(\rho_\e) +(1-\th) \Phi(\rho), \cd\Big)- \frac{\delta U}{\delta\nu} \big( \Phi(\rho), \cd\big)\|_n  d\th \to 0, ~ \mbox{as}~ \e\to 0,
\eeaa
thanks to the continuity of ${\d U\over \d\nu}$ in $\nu$ under $\|\cd\|_n$. This, together with \reff{I1e}, leads to
\beaa
\lim_{\e\to 0}{1\over \e}\big[ \wh U(\rho_\e) - \wh U(\rho) \big] = \int_{\dbR^k}\int_{\dbR^k} \frac{\delta U}{\delta\nu} \Big(\Phi(\rho), y\Big) {\d \Phi\over \d\rho}(\rho, x; dy) (\rho'-\rho)(dx),
\eeaa
which implies \reff{linearchain} immediately.
\qed

\begin{rem}
\label{rem-derivative}
By considering generalized derivatives in appropriate dual space, we may define higher order derivatives of $\Phi$, including the Lions derivative $\pa_\rho \Phi(\rho, x):= \pa_x {\d \Phi\over \d \rho}(\rho, x)$. Alternatively, since later on we will always consider certain composite  function $\wh U$, we may define higher order derivatives through the left side of \reff{linearchain}. 
\end{rem}
 
 \section{The master equation and the monotonicities}
\label{sect-master} 
\setcounter{equation}{0}
Throughout the paper, Assumption \ref{assum-fix} will always be in force. Denote
\bea
\label{barH}
\wh H(x, p, \rho) := H(x, p, \Phi(\rho)),\q (x, p, \rho)\in \dbR^d\times \dbR^d \times \cP_2(\dbR^{2d}).
\eea
The derivatives of $\wh H$ with respect to $\rho$ are understood as in Proposition \ref{prop-chain} and Remark \ref{rem-derivative}. Then \reff{X*} becomes
\bea
\label{barX*}
\left.\ba{c}
\dis X_t^{*}=\xi+\int_{t_0}^{t} \pa_p \wh H\Big(X_s^{*}, \pa_x V(s, X^*_s, \mu^*_s), \rho^*_s\Big)ds+B_t^{t_0}+\beta B_t^{0,t_0},\\
\mbox{where}\q \mu^*_s:= \cL_{X^*_s|\cF^0_s},\q \rho^*_s:= \cL_{(X^*_s, \pa_x V(s, X^*_s, \mu^*_s))|\cF^0_s}.
\ea\right.
\eea
On the other hand, it follows from the standard stochastic control theory that, for given $t_0, \mu$, the optimization \reff{MFE} is associated with the following Backward SDE: recalling \reff{Hp},
\bea
\label{Y*}
\left.\ba{c}
\dis Y_t^{*}= G(X^*_T, \mu_T^*) - \int_t^T Z^*_s dB_s -\int_t^T Z^{0,*}_s dB_s^0\\
\dis  +\int_{t}^T \Big[\wh H(\cd) - \pa_x V(s, X^*_s, \mu^*_s) \cd \pa_p \wh H(\cd)\Big]\Big(X_s^{*}, \pa_x V(s, X^*_s,\mu^*_s), \rho^*_s\Big)ds,
\ea\right.
\eea
which, together with \reff{barX*}  form the MFGC system. We note that this is the SDE counterpart of the MFGC system \reff{SPDEint}.  In particular, we have
\bea
\label{YVX}
Y^*_t = V\big(t, X^*_t, \mu_t^*\big).
\eea
Then, by applying the It\^{o}'s formula (c.f. \cite[Theorem 4.17]{CD2},\cite{BLPR,CCD}) we obtain
\bea
&&d V(t, X^*_t, \mu^*_t) =  \Big[\pa_t V + \pa_x V\cdot \pa_p \wh H(X^*_t, \pa_x V, \rho^*_t) + \frac{1+\b^2}{2} \tr\big(\pa_{xx} V \big)\Big](t, X^*_t, \mu^*_t) dt \nonumber\\
&&+\pa_xV(t,X^*_t,\mu^*_t)\cd dB_t  + \beta\Big[\pa_xV(t,X^*_t,\mu^*_t) + \tilde{\mathbb E}_{\cF_t}\big[\pa_\mu V(t,X^*_t,\mu^*_t,\tilde X^*_t) \big]\Big]\cd dB_t^0
\nonumber\\
&&+\tr \Big(\tilde \dbE_{\cF_t}\big[\pa_\mu V(t,X^*_t,\mu_t^*,\tilde X_t^*) (\pa_p\wh H(t,\tilde X_t^*, \pa_x V(t, \tilde X_t, \mu_t^*), \rho^*_t))^\top\big]\Big) dt  \label{Ito}\\
&&+ \tr \Big(\b^2\tilde \dbE_{\cF_t}\big[\pa_x\pa_\mu V(t,X^*_t,\mu^*_t,\tilde X^*_t)+\frac{1+\beta^2}{2} \pa_{\tilde x}\pa_\mu V(t, X^*_t, \mu_t^*, \tilde X^*_t)\big]\nonumber\\
&&\qq +\frac{\beta^2}{2}\bar{\tilde \dbE}_{\cF_t}\big[\pa_{\mu\mu}V(t,X^*_t,\mu^*_t,\tilde X^*_t,\bar X^*_t) \big]\Big)dt.\nonumber
\eea
Here as usual $\tilde X^*,\bar X^*$ are conditionally independent copies of $X^*$, conditional on $\dbF^0$. 
 Comparing this with \eqref{Y*}, we derive the master equation:  for independent copies $\xi, \tilde \xi, \bar \xi$ with law $\mu$,
\bea
\label{master}
\left.\ba{c}
\dis \cL V(t,x,\mu) :=\pa_t V + \frac{\h\b^2}{2} \tr(\pa_{xx} V) + \wh H(x,\partial_x V,\mathcal{L}_{(\xi,\pa_xV(t,\xi,\mu))}) + \cM V =0, \\
\dis V(T,x,\mu) = G(x,\mu),\q\mbox{where}\\
\dis \cM V(t,x,\mu)  := \tr\Big( \bar{\tilde \dbE}\Big[\frac{\h\b^2}{2} \pa_{\tilde x} \pa_\mu V(t,x, \mu, \tilde \xi) +\b^2\pa_x\pa_\mu V(t,x,\mu,\tilde \xi)+\frac{\b^2}{2}\pa_{\mu\mu}V(t,x,\mu,\bar\xi,\tilde\xi) \\
\dis   + \pa_\mu V(t, x, \mu, \tilde \xi)(\pa_p\wh H)^\top(\tilde \xi,\pa_x V(t, \tilde \xi, \mu), \mathcal{L}_{(\xi,\pa_xV(t,\xi,\mu))})
\Big]\Big), \q\mbox{and}\q \h \b^2 := 1+\b^2.
\ea\right.
\eea

In addition to Assumption \ref{assum-fix}, we assume
\begin{assum} \label{assum-regH} 
$\wh H\in \cC^2(\mathbb R^{2d}\times\mathcal{P}_2(\mathbb R^{2d}))$ with bounded $\pa_{xp}\wh H, \pa_{xx}\wh H, \pa_{pp}\wh H, \pa_{x\mu}\wh H, \pa_{p\mu}\wh H$.
 \end{assum}
 Since we will work on the master equation, here we impose our conditions directly on $\wh H$, rather than on $b, f$. It is not hard to find some sufficient conditions on $b$ and $f$ to ensure these.

\subsection{The monotonicities}

In this subsection we introduce three types of monotonicity conditions: Lasry-Lions monotonicity, displacement semi-monotonicity, and anti-monotonicity.   
\begin{defn}
Assume $U\in\mathcal{C}^1(\mathbb R^d\times\mathcal{P}_2(\mathbb R^d))$ and $\pa_\mu U(\cd,\mu,\tilde x)\in \cC^1(\mathbb R^d)$ for all $(\mu,\tilde x)\in\mathcal{P}_2(\mathbb R^d)\times\mathbb R^d$. We say $U$ is Lasry-Lions monotone if
\begin{equation}
\label{LL}
MON^{LL} U(\xi, \eta):=\tilde \dbE\Big[\big\langle \pa_{x\mu }U(\xi,\mathcal{L}_\xi,\tilde \xi)\tilde\eta,\eta\big\rangle\Big]\ge 0,\q \forall \xi,\eta\in \mathbb L^2(\cF^1_T).
\end{equation}
\end{defn}
We note that, since $(\xi,\eta)$ is $\cF^1_T$-measurable, so here $(\tilde \xi, \tilde \eta)$ is an independent copy (instead of conditionally independent copy).

\begin{defn}
Assume $U,\pa_x U\in\mathcal{C}^1(\mathbb R^d\times\mathcal{P}_2(\mathbb R^d))$. For any $\l\ge 0$, we say $U$ is displacement $\l$-monotone if, for all $\xi,\eta\in \mathbb L^2(\cF^1_T)$
\begin{equation}
\label{lambdadisplacementsm}
MON^{disp}_\l U(\xi, \eta):=\tilde\dbE\Big[\big\langle \pa_{x\mu }U(\xi,\mathcal{L}_\xi,\tilde \xi)\tilde\eta,\eta\big\rangle+\big\langle\pa_{xx}U(\xi,\mathcal{L}_{\xi})\eta,\eta\big\rangle + \l |\eta|^2\Big] \ge 0.
\end{equation}
In particular, we say $U$ is displacement monotone when $\l=0$, and displacement semi-monotone if it is displacement $\l$-monotone for some $\l>0$.
 \end{defn}

Moreover, denote 
\bea
\label{D4}
D_4 := \Big\{\vec \l = (\l_0, \l_1, \l_2, \l_3): \l_0>0, \l_1\in \dbR, \l_2 >0, \l_3\ge 0\Big\}.
\eea

\begin{defn}
\label{defn-anti}
Let $U\in \cC^2(\dbR^d\times \cP_2(\mathbb R^d))$ and $\vec\l\in D_4$. We say $U$ is \mbox{$\vec\l$-anti-monotone} if,
\bea
\label{anti}
\left.\ba{c}
MON^{anti}_{\vec \l} U(\xi, \eta):= \tilde\dbE\bigg[\l_0\langle\pa_{xx}U(\xi,\mathcal{L}_{\xi})\eta,\eta\rangle+\l_1\langle\pa_{x\mu}U(\xi,\mathcal{L}_{\xi},\tilde\xi)\tilde\eta,\eta\rangle+\left|\pa_{xx}U(\xi,\mathcal{L}_{\xi})\eta\right|^2\\
+\l_2\big|\tilde{\mathbb E}_\cF[\pa_{x\mu}U(\xi,\mathcal{L}_{\xi},\tilde\xi)\tilde \eta]\big|^2 - \l_3|\eta|^2 \bigg]\leq 0,\q \forall \xi,\eta\in \mathbb L^2(\cF^1_T).
\ea\right.
\eea
\end{defn}

\begin{rem}
\label{rem-mon}
(i) By \cite[Remark 2.4]{GMMZ}, Lasry-Lions monotonicity and displacement monotonicity are equivalent to the following forms, respectively, which are more often seen in the literature:
\beaa
&\dis \dbE\Big[ U(\xi_1, \cL_{\xi_1}) + U(\xi_2, \cL_{\xi_2})  - U(\xi_1, \cL_{\xi_2}) -U(\xi_2, \cL_{\xi_1}) \Big] \ge 0,\q \forall \xi_1, \xi_2\in \dbL^2(\cF^1_T);\\
&\dis \dbE\Big[ \big\langle \pa_xU(\xi_1, \cL_{\xi_1}) -\pa_x U(\xi_2, \cL_{\xi_2}),~ \xi_1-\xi_2\big\rangle  \Big] \ge 0, \q \forall \xi_1, \xi_2\in \dbL^2(\cF^1_T).
\eeaa

(ii) Consider the case that $\pa_x U(x, \mu) = \pa_\mu \cU(\mu, x)$ for some  $\mathcal{U}\in\cC^2(\cP_2(\dbR^d))$. Then the Lasry-Lions monotonicity of $U$ is equivalent to the convexity of the mapping $\mu \in \cP_2(\dbR^d) \mapsto \cU(\mu)$, and  the displacement monotonicity of $U$ is equivalent to the convexity of  the mapping $\xi \in \dbL^2(\cF^1_T) \mapsto \cU(\cL_\xi)$, see e.g. \cite{CDLL,CD1}.

(iii) Both the Lasry-Lions monotonicity  (provided $\pa_{xx} U$ is bounded) and the displacement monotonicity imply the displacement semi-monotonicity. However, the Lasry-Lions monotonicity  and the displacement monotonicity do not imply each other, see \cite[Remark 2.5]{GMMZ}.

(iv) By setting $\l_0=\l_1=\l_2=1$ and $\l_3=0$, \reff{anti} implies 
\beaa
\tilde\dbE\bigg[\langle\pa_{xx}U(\xi,\mathcal{L}_{\xi})\eta,\eta\rangle+\langle\pa_{x\mu}U(\xi,\mathcal{L}_{\xi},\tilde\xi)\tilde\eta,\eta\rangle \bigg]\leq 0,
\eeaa
which is in the opposite direction of \reff{lambdadisplacementsm} with $\l=0$. Moreover, if $\pa_{xx}U$ is non-negative definite, then we further have
\beaa
\tilde\dbE\bigg[\langle\pa_{x\mu}U(\xi,\mathcal{L}_{\xi},\tilde\xi)\tilde\eta,\eta\rangle \bigg]\leq 0,
\eeaa
which is  in the opposite direction of \eqref{LL}. That's why we call \reff{anti} anti-monotonicity.

(v) If $U$ satisfies \eqref{lambdadisplacementsm} for some $\l\geq 0$, then $\pa_{xx}U+\l I$ is non-negative definite, see \cite[Lemma 2.6]{GMMZ}.
\end{rem}

\subsection{A road map towards the global wellposedness}
\label{sect-strategy}
Our ultimate goal is to establish the global wellposedness of the master equation \reff{master}. We shall adopt the strategy in \cite{GMMZ, MZ3}, which consists of three steps:

\vskip 2pt
{\it Step 1.} Introduce appropriate monotonicity condition on data which ensures the propagation of the monotonicity, one of the three types introduced in the previous subsection, along any classical solution to the master equation. 

\vskip 2pt
{\it Step 2.} Show that the monotonicity of $V(t,\cdot,\cdot)$ implies an (a priori) uniform Lipschitz continuity of $V$ in the measure variable $\mu$.

\vskip 2pt
{\it Step 3.} Combine the local wellposedness of classical solutions and the above uniform Lipschitz continuity to obtain the global wellposedness of classical solutions.

Moreover, following \cite{CDLL} one may continue to investigate the convergence problem:

\vskip 2pt
{\it Step 4.} Use the classical solution $V$ to prove the convergence of the related $N$-player game. 

\vskip 2pt
In this paper we shall focus on  {\it Step 1} only, and we leave the remaining  three steps to future research. We emphasize that {\it Step 1} (and {\it Step 2}) considers a prior estimates, and thus throughout the paper we shall also assume:

\begin{assum} \label{assum-regV} 
$V\in \mathcal{C}^{1,2,2}([0,T]\times\mathbb R^d\times\mathcal{P}_2(\mathbb R^d))$ is a classical solution of the master equation \reff{master} such that $\pa_{xx}V(t,\cd,\cd)\in\mathcal{C}^2(\dbR^d\times\mathcal{P}_2(\mathbb R^d))$, $\pa_{x\mu}V(t,\cd,\cd,\cd)\in\mathcal{C}^2(\dbR^d\times\mathcal{P}_2(\mathbb R^d)\times \dbR^d)$,  and all the second and higher order derivatives of $V$ involved above are uniformly bounded and continuous in $t$.
 \end{assum}
 We note that we do not require $V$ or its first order derivatives to be bounded. Moreover,  since $G = V(T,\cd,\cd)$, so the above assumption also ensures the regularity of $G$. We shall also remark that, Assumption \ref{assum-regV} is about the existence of classical solutions of the master equation \reff{master}, which implies the uniqueness of the mean field equilibrium (c.f. \cite[Remark 2.10 (ii)]{GMMZ}). The uniqueness of classical solutions satisfying the desired Lipschitz continuity is standard, see e.g. the arguments in \cite[Theorem 6.3]{GMMZ}.

\section{Propagation of Lasry-Lions monotonicity}
\label{sect-LL}
\setcounter{equation}{0}

To propagate the Lasry-Lions monotonicity of $V$, we impose the following assumption on $\wh H$.
\begin{assum}\label{assum-LL} 
For any $\xi, \eta, \g,  \zeta\in\mathbb L^2(\cF^1_T)$ and $\varphi: \dbR^d\to \dbR^d$ Lipschitz continuous,
\bea
\label{HLL}
\left.\ba{c}
\dis \tilde\dbE\bigg[    \Big\langle \zeta,  \wh H_{pp}(\xi)\zeta\Big\rangle - \Big\langle \eta, ~\wh H_{x\rho_1}(\xi,\tilde \xi)  \tilde \eta+ \wh H_{x\rho_2}(\xi,\tilde \xi)[ \tilde\g +\tilde \zeta]\Big\rangle \\
\dis\qq -\Big\langle \g -\zeta,  ~\wh H_{p\rho_1}(\xi, \tilde \xi) \tilde \eta+\wh H_{p\rho_2}(X_t,\tilde X_t)[ \tilde \g+\tilde \zeta] \Big\rangle  \bigg] \le 0,
\ea\right.
\eea
where $\wh H_{pp}(x) := \pa_{pp} \wh H\Big(x, \f(x), \cL_{(\xi, \f(\xi))}\Big)$, $\wh H_{x\rho}(x,\tilde x) := \pa_{x\rho}\wh H\Big(x, \f(x), \cL_{(\xi, \f(\xi))}, \tilde x, \f(\tilde x)\Big)$, and similarly for $\wh H_{p\rho}(x,\tilde x)$.

\end{assum}

Our main result of this section is:
\begin{thm}\label{thm:LL}
Let Assumptions \ref{assum-fix},  \ref{assum-regH}, \ref{assum-regV}, and \ref{assum-LL} hold. If $G$ satisfies the Lasry-Lions monotonicity \reff{LL}, then $V(t,\cd,\cd)$ satisfies \reff{LL} for all $t\in[0,T]$.
\end{thm}
\proof 
Without loss of generality, we shall prove the theorem only for $t=0$.

For any $\xi,\eta\in\mathbb L^2(\mathcal{F}_0)$, inspired by \reff{barX*} we consider the following system of McKean-Vlasov SDEs, which clearly has a unique solution $(X,\delta X)$ under Assumptions \ref{assum-regH} and \ref{assum-regV}:
\bea
\label{XY}
\left.\ba{lll}
\dis X_t = \xi +\int_0^t \wh H_p(X_s, \pa_x V(s, X_s, \mu_s), \rho_s) ds +  B_t+\beta B_t^0;\\
\dis \delta X_t = \eta +\int_0^t \Big[\wh H_{px} (X_s) \delta  X_s + \wh H_{pp}(X_s)[\Gamma_s+\Upsilon_s ]  + N_s\Big]ds;\q\mbox{where}\\
\dis \mu_t:=\mathcal{L}_{X_t|\mathcal{F}_t^{0}},\q \rho_t := \cL_{(X_t,\pa_xV(t,X_t,\mu_t))|\mathcal{F}_t^{0}};\ms\\
\dis \Gamma_t:=\pa_{xx}V(X_t)\delta X_t,\q \Upsilon_t:=\tilde \dbE_{\mathcal{F}_t}[\pa_{x\mu}V(X_t,\tilde X_t)\delta \tilde X_t];\ms\\
\dis N_t :=\tilde \dbE_{\mathcal{F}_t}\Big[ \wh H_{p\rho_1}(X_t,\tilde X_t) \delta  \tilde X_t+ \wh H_{p\rho_2}(X_t,\tilde X_t) \big[\tilde\Gamma_t+\tilde \Upsilon_t\big]\Big].
\ea\right.
\eea
Here $(\tilde X, \d\tilde X, \tilde \G, \tilde \Upsilon)$ is a conditionally independent copy of $(X, \d X, \G, \Upsilon)$, conditional on $\dbF^0$.
Moreover, here and in the sequel,  for simplicity of notation, we omit the variables $(t, \mu_t)$ inside $V$ and its derivatives, and omit $\rho_t$ and  $\partial_x V$ inside $\wh H$ and its derivatives, for example, 
\bea
\label{omit}
\left.\ba{c}
\pa_{x\mu} V(X_t, \tilde X_t) = \pa_{x\mu}V(t, X_t, \mu_t, \tilde X_t),\q \wh H_p(X_t) := \pa_p \wh H(X_t, \pa_x V(t, X_t, \mu_t) ,\rho_t ), \\ 
\wh H_{p\rho}(X_t, \tilde X_t):= (\wh H_{p\rho_1},\wh H_{p\rho_2})(X_t, \tilde X_t):= \pa_{p\rho}\wh H(X_t,\pa_x V(t, X_t, \mu_t), \rho_t,  \tilde X_t, \pa_x V(t, \tilde X_t, \mu_t)).
\ea\right.
\eea

Introduce,
\bea
\label{It}
I(t):=\mathbb E[\langle\Upsilon_t,\delta X_t\rangle] = MON^{LL} V(t,\cd,\cd)(X_t, \d X_t).
\eea
Apply It\^{o} formula \reff{Ito} and since $V$ satisfies the master equation \reff{master}, we get
\bea\label{dtI-LL}
\left.\ba{c}
\dis{d\over dt} I(t)=\tilde\dbE\bigg[    \Big\langle \Upsilon_t,  ~\wh H_{pp}(X_t)\Upsilon_t\Big\rangle - \Big\langle \delta  X_t, ~\wh H_{x\rho_1}(X_t,\tilde X_t) \delta \tilde X_t+ \wh H_{x\rho_2}(X_t,\tilde X_t)[\tilde \Gamma_t+\tilde \Upsilon_t]\Big\rangle \\
\dis -\Big\langle \Gamma_t- \Upsilon_t,  ~\wh H_{p\rho_1}(X_t, \tilde X_t)\delta \tilde X_t+\wh H_{p\rho_2}(X_t,\tilde X_t)\big[\tilde\Gamma_t+\tilde \Upsilon_t] \Big\rangle  \bigg].
\ea\right.
\eea
The calculation is lengthy but quite straightforward, we postpone the details to Appendix. 
Take conditional expectation on $\cF^0_t$, then by the desired conditional independence we may apply  \eqref{HLL} to obtain:
\bea
\label{dI}
{d\over dt} I(t) \le 0.
\eea
Note that, by the Lasry-Lions monotonicity of $G = V(T,\cd,\cd)$, we have $I(T)\ge 0$. Then \reff{dI} clearly implies $I(0) \ge 0$, and hence $V(0,\cd,\cd)$ satisfies the Lasry-Lions monotonicity \reff{LL}.
\qed

\begin{rem}
In \eqref{XY} $X$ is the agent's state process along the (unique) mean field equilibrium, and $\d X$ is the gradient of $X$ when its initial condition $\xi$ is perturbed along the direction $\eta$.
\end{rem}

\begin{rem}
\label{rem-LL}
Note that \reff{dtI-LL} is an equality, so our condition \reff{HLL} is essentially sharp for the propagation of Lasry-Lions monotonicity, in particular for \reff{dI}. In \cite{CardaliaguetLehalle,CD1, K0} the uniqueness of the mean field game system is obtained when $b(\cdot,a,\cdot) = a$ (or slightly more general form), and $f$ satisfies the Lasry-Lions monotonicity in the following sense: for any $\xi_i, \a_i\in \cL^2(\cF)$, $i=1,2$, 
\bea
\label{fmon}
\mathbb E\big[f(\xi^1,\alpha^1,\mathcal{L}_{(\xi^1,\alpha^1)})+f(\xi^2,\alpha^2,\mathcal{L}_{(\xi^2,\alpha^2)})-f(\xi^1,\alpha^1,\mathcal{L}_{(\xi^2,\alpha^2)})-f(\xi^2,\alpha^2,\mathcal{L}_{(\xi^1,\alpha^1)})\big]\geq 0.
\eea
We claim that in this case \reff{dI} holds true, and hence the Lasry-Lions monotonicity propagates. We postpone its proof to the Appendix.
\end{rem}

\begin{rem}
\label{rem-LL-separable}
For the standard MFG with $b(x,a,\nu) = a$ (and $f= f(x, a, \mu)$), it is observed in \cite{GMMZ} that it is hard to propagate the Lasry-Lions monotonicity unless $f$ is separable: $f(x, a, \mu) = f_0(x, a) + f_1(x, \mu)$. The dependence on the law of $\a$ in MFGC actually helps for the propagation of the Lasry-Lions monotonicity. In particular, in this case we do not require $f$ to be separable.
\end{rem}

We next provide an example with a more general $b$, which does not seem to be covered by the analysis of mean field game systems (or master equations) in the literature.

\begin{eg}
\label{eg-LasryLions}
We consider a special case of \reff{nonseparable} with $d=1$: 
\bea
\label{nonseparable2}
\left.\ba{c}
\dis b(x, a, \cL_{(\xi,\a)}) = -a  + b_1(\dbE[\xi],\dbE[\a]) + b_2(x), \ms\\
\dis f(x, a, \cL_{(\xi,\a)}) = {|a|^2\over 2} -  c_1 a \dbE[\a]  + c_2 x\dbE[\xi] + c_3 x \dbE[\a]  + f_1(x),
\ea\right.
\eea
where $0<c_1<1$ and $c_2, c_3>0$ are constants. Assume the matrix
\bea
\label{matrix1}
\begin{bmatrix}
1-[\bar c_1\pa_{m_2} b_1 -\hat c_1]& 0& {1\over 2}[\hat c_3 -  \pa_{m_1} b_1 ] \\
0 & [\bar c_1\pa_{m_2} b_1 -\hat c_1]  &  {1\over 2}[\hat c_3 +  \pa_{m_1} b_1 ]\\
 {1\over 2}[\hat c_3 -  \pa_{m_1} b_1 ]&   {1\over 2}[\hat c_3 +  \pa_{m_1} b_1 ]& c_2
\end{bmatrix} \ge 0,
\eea
where $ \hat c_1 := {c_1\over 1-c_1}$, $\bar c_1:=\frac{1}{1-c_1}$,  $\hat c_3 := {c_3\over 1-c_1}$, and $m_1, m_2$ stand for $\dbE[\xi], \dbE[\a]$. Then \reff{HLL} holds true.
\end{eg}
\proof By Example \ref{eg-nonseparable} we see that
\beaa
&\dis \Phi(\cL_{(\xi, \eta)}) = \cL_{(\xi, ~ \hat c_1\dbE[\eta]  + \eta)},\\
&\dis H(x, p, \cL_{(\xi,\a)}) = -{1\over 2} \big|c_1\dbE[\a]  + p\big|^2 + p\big[b_1(\dbE[\xi], \dbE[\a]) + b_2(x)\big]   +c_2 x\dbE[\xi] +  c_3 x \dbE[\a]  + f_1(x).
\eeaa
Note that $\dbE[\a] = [1+\hat c_1] \dbE[\eta] = \bar c_1 \dbE[\eta]$. Then
\beaa
 \wh H(x, p, \cL_{(\xi, \eta)}) =  -{1\over 2} \big|\hat c_1\dbE[\eta]   +p\big|^2 +p\big[ b_1(\dbE[\xi], \bar c_1\dbE[\eta])+b_2(x)\big] +c_2 x\dbE[\xi]+  \hat c_3 x \dbE[\eta]  + f_1(x).
\eeaa
One may compute straightforwardly that
\bea
\label{whHpp}
\wh H_{pp} = -1,\q  \wh H_{x\rho_1} = c_2,\q \wh H_{x\rho_2} =  \hat c_3,\q \wh H_{p \rho_1}= \pa_{m_1} b_1,\q \wh H_{p\rho_2}=\bar c_1\pa_{m_2}b_1-\hat c_1.
\eea
Then, noting that $ \pa_{m_1} b_1$ and $ \pa_{m_1} b_2$ are deterministic, 
\beaa
\left.\ba{lll}
\dis  \tilde\dbE\Big[   - \wh H_{pp}(\xi)|\zeta|^2 + \eta\big[\wh H_{x\rho_1}(\xi,\tilde \xi)  \tilde \eta+ \wh H_{x\rho_2}(\xi,\tilde \xi)[ \tilde\g +\tilde \zeta]\big] \\
\dis\qq +[\g -\zeta]\big[\wh H_{p\rho_1}(\xi, \tilde \xi) \tilde \eta+\wh H_{p\rho_2}(X_t,\tilde X_t)[ \tilde \g+\tilde \zeta] \big]  \Big] \\
\dis=\tilde\dbE\Big[|\zeta|^2 + c_2 \eta \tilde \eta  + \hat c_3 \eta [ \tilde\g +\tilde \zeta]\big]  +[\g -\zeta]\big[\pa_{m_1} b_1\tilde \eta+ [\bar c_1\pa_{m_2} b_1 -\hat c_1][ \tilde \g+\tilde \zeta] \big]  \Big]\ms\\
\dis = \dbE[|\zeta|^2] + c_2 \big|\dbE[\eta]\big|^2 + \hat c_3 \dbE[\eta] \big[\dbE[\g] + \dbE[\zeta]\big] + \pa_{m_1} b_1 \dbE[\eta]\big[\dbE[\g]-\dbE[\zeta]\big]\ms \\
\dis\qq +[\bar c_1\pa_{m_2} b_1 -\hat c_1] \big[ |\dbE[\g]|^2 - |\dbE[\zeta]|^2\big]\ms\\
\dis \ge \big[1- [\bar c_1\pa_{m_2} b_1 -\hat c_1]\big] \big|\dbE[\zeta]\big|^2 + [\bar c_1\pa_{m_2} b_1 -\hat c_1]  \big|\dbE[\g]\big|^2 +  c_2 \big|\dbE[\eta]\big|^2\ms\\
\dis \q + \big[\hat c_3 +  \pa_{m_1} b_1 \big] \dbE[\eta] \dbE[\g]  +  \big[\hat c_3 -  \pa_{m_1} b_1 \big] \dbE[\eta] \dbE[\zeta].
\ea\right.
\eeaa
This, together with  \reff{matrix1}, clearly implies  \reff{HLL}.
\qed

\section{Propagation of displacement $\l$-monotonicity}
\label{sect-displacement}
\setcounter{equation}{0}

In this section we fix a constant $\l\ge 0$.

\begin{assum}
\label{assum-displacement} 
For any $\xi, \eta, \g,  \zeta\in\mathbb L^2(\cF^1_T)$ and $\varphi: \dbR^d\to \dbR^d$ Lipschitz continuous,
\bea
\label{Hdisplacement}
\left.\ba{lll}
\dis \tilde\dbE\bigg[    \Big\langle\g + \zeta,  \wh H_{pp}(\xi)[\g + \zeta]\Big\rangle -  \Big\langle \eta,  [\wh H_{xx}(\xi)-2\l \wh H_{px}(\xi)]\eta\Big\rangle \\
\dis+\Big\langle \g + \zeta,  ~\big[\wh H_{p\rho_1}(\xi, \tilde \xi)+\wh H_{\rho_2 x}(\tilde \xi,\xi)+2\l \wh H_{\rho_2 p}(\tilde \xi,\xi)\big]\tilde\eta+2\l \wh H_{pp}(\xi)\eta\Big\rangle\\
\dis +\Big\langle\g + \zeta, \wh H_{p\rho_2}(\xi,\tilde \xi)[\tilde \g + \tilde \zeta]\Big]\Big\rangle- \Big\langle \eta,  ~[\wh H_{x\rho_1}(\xi,\tilde \xi)-2\l \wh H_{p\rho_1}(\xi,\tilde \xi)] \tilde \eta\Big\rangle\bigg]\leq 0,
\ea\right.
\eea
where $\wh H_{pp}, \wh H_{x\rho}, \wh H_{p\rho}$ are as in Assumption \ref{assum-LL}.
\end{assum}

\begin{thm}\label{thm:displacement}
Let Assumptions \ref{assum-fix},  \ref{assum-regH}, \ref{assum-regV}, and \ref{assum-displacement}  hold. If $G$ satisfies the displacement $\l$-monotonicity \reff{lambdadisplacementsm}, then $V(t,\cd,\cd)$ satisfies \reff{lambdadisplacementsm} for all $t\in[0,T]$.
\end{thm}
\proof 
Without loss of generality, we shall prove the theorem only for $t=0$. We will continue to use the notation in the proof of Theorem \ref{thm:LL}.

Introduce
\begin{equation}\label{barI}
\bar I(t):=\mathbb E[\langle\Gamma_t,\delta X_t\rangle],\q\mbox{and thus}\q I(t)+\bar I(t) + \l\mathbb E[|\delta X_t|^2] = MON^{disp}_\l V(t,\cd,\cd)(X_t, \d X_t).
\end{equation}
Similarly to \reff{dtI-LL} we can show that, again see more details in Appendix, 
\bea\label{cdbarI}
\left.\ba{c}
\dis {d\over dt} {\bar I}(t)  =\tilde\dbE\bigg[  \Big\langle \wh H_{pp}(X_t)\Gamma_t, \Gamma_t\Big\rangle+ 2 \Big\langle \wh H_{pp}(X_t) \Gamma_t, \Upsilon_t\Big]\Big\rangle\\
\dis+ 2\Big\langle\Gamma_t, \wh H_{p\rho_1}(X_t, \tilde X_t)\delta  \tilde X_t+\wh H_{p\rho_2}(X_t,\tilde X_t)[\tilde \Gamma_t+\tilde \Upsilon_t]\Big\rangle - \Big\langle \wh H_{xx}(X_t) \delta  X_t, \delta  X_t\Big\rangle\bigg].
\ea\right.
\eea
Moreover,  by  \eqref{XY} we have
\begin{equation}\label{deltaX2}
\frac{d}{dt}\mathbb E\left[|\delta X_t|^2\right]=2\mathbb E\Big[\Big\langle \wh H_{px}(X_t)\delta X_t+\wh H_{pp}(X_t)[ \Upsilon_t+\Gamma_t]+N_t,~\delta X_t \Big\rangle\Big].
\end{equation}
Combining \eqref{dtI-LL}, \eqref{cdbarI}, and \eqref{deltaX2}, and recalling the $N$ in \reff{XY}, we deduce that
\bea
\label{cdI}
\left.\ba{lll}
\dis {d\over dt} \Big[MON^{disp}_\l V(t,\cd,\cd)(X_t, \d X_t)\Big] = {d\over dt}\Big[I(t)+\bar I(t) + \l\mathbb E[|\delta X_t|^2]\Big]\\
\dis =\tilde\dbE\bigg[    \Big\langle \Upsilon_t+\Gamma_t,  \wh H_{pp}(X_t)[\Upsilon_t+\Gamma_t]+\wh H_{p\rho_2}(X_t,\tilde X_t)[\tilde \Gamma_t+\tilde \Upsilon_t]\Big\rangle  \\
\dis\q +\Big\langle  \Upsilon_t+ \Gamma_t,  [\wh H_{p\rho_1}(X_t, \tilde X_t)+\wh H_{\rho_2x}(\tilde X_t,X_t)+2\l \wh H_{\rho_2p}(\tilde X_t,X_t)]\delta \tilde X_t+2\l \wh H_{pp}(X_t)\delta X_t\Big\rangle\\
\dis \q- \Big\langle \delta  X_t, [\wh H_{x\rho_1}(X_t,\tilde X_t)-2\l \wh H_{p\rho_1}(X_t,\tilde X_t)] \delta \tilde X_t+[\wh H_{xx}(X_t)-2\l \wh H_{px}(X_t)]\delta X_t\Big\rangle\bigg].
\ea\right.
\eea
Then, by the desired conditional independence of the involved processes above, conditional on $\cF^0_t$, we have
\bea
\label{dIbar}
{d\over dt} \Big[MON^{disp}_\l V(t,\cd,\cd)(X_t, \d X_t)\Big] \le 0.
\eea
Note that $V(T,\cd, \cd)=G$ satisfies \reff{lambdadisplacementsm}, then clearly $V(0,\cd,\cd)$ also satisfies \reff{lambdadisplacementsm}.
\qed

\ms

We next provide a sufficient condition for  Assumption \ref{assum-displacement}. Denote, for any $A\in \dbR^{d\times d}$,
\bea\label{kappaA}
\left.\ba{c}
\dis  |A| :=  \sup_{|x|=|y|=1} \langle Ax, y\rangle,\q \underline \k(A) := \inf_{|x|=1} \langle Ax, x\rangle = \mbox{the smallest eigenvalue of ${1\over 2}[A+A^\top]$},\\
\dis  \ol \k(A):=\sup_{|x|=1} \langle Ax, x\rangle=-\ul\k(-A).
\ea\right.
\eea

\begin{prop}
\label{prop-dissufficient}
Assume there exists a constant $c_0\ge 0$ such that $|\pa_{p\rho_2}\wh H|\le c_0$, and $\wh H_{pp} < - c_0 I_d$, where $I_d$ denotes the $d\times d$ identity matrix. Then the following condition implies \reff{Hdisplacement}:
\bea
\label{lambdaHH}
\left.\ba{c}
\dis\tilde\dbE\bigg[  \Big\langle \eta,  [\wh H_{xx}(\xi)-2\l \wh H_{px}(\xi)]\eta\Big\rangle  + \Big\langle \eta,  ~[\wh H_{x\rho_1}(\xi,\tilde \xi)-2\l \wh H_{p\rho_1}(\xi,\tilde \xi)] \tilde \eta\Big\rangle - {|\L(\xi, \eta)|^2\over 4}\bigg] \ge 0, \\
\dis \mbox{where}\q \L(\xi, \eta) := (- \wh H_{pp}(\xi) - c_0I_d)^{-{1\over 2}} \times\ms\\
\dis \Big[ \tilde \dbE_{\cF^1_T}\big[[\wh H_{p\rho_1}(\xi, \tilde \xi)+\wh H_{\rho_2 x}(\tilde \xi,\xi)+2\l \wh H_{\rho_2 p}(\tilde \xi,\xi)]\tilde\eta\big]+2\l \wh H_{pp}(\xi)\eta\Big],
\ea\right.
\eea
for all $\xi, \eta\in \dbL^2(\cF_T^1)$.
In particular, when $\l=0$, the above reduces to:
 \bea
\label{lambdaHH2}
\left.\ba{c}
\dis\tilde\dbE\bigg[  \Big\langle \eta,  \wh H_{xx}(\xi)\eta\Big\rangle  + \Big\langle \eta,  \wh H_{x\rho_1}(\xi,\tilde \xi) \tilde \eta\Big\rangle - {|\L(\xi, \eta)|^2\over 4}\bigg] \ge 0, \ms\\
\dis \mbox{where}\q \L(\xi, \eta) := (- \wh H_{pp}(\xi) - c_0I_d)^{-{1\over 2}} \tilde \dbE_{\cF_T^1}\Big[[\wh H_{p\rho_1}(\xi, \tilde \xi)+\wh H_{\rho_2 x}(\tilde \xi,\xi)]\tilde\eta\Big].
\ea\right.
\eea
\end{prop}
\proof Denote $\Xi:= \big[\wh H_{p\rho_1}(\xi, \tilde \xi)+\wh H_{\rho_2 x}(\tilde \xi,\xi)+2\l \wh H_{\rho_2 p}(\tilde \xi,\xi)\big]\tilde\eta$. Note that
\beaa
&&\dis \tilde\dbE\bigg[    \Big\langle\g + \zeta,  \wh H_{pp}(\xi)[\g + \zeta]\Big\rangle+\Big\langle \g + \zeta,  \Xi+2\l \wh H_{pp}(\xi)\eta\Big\rangle+\Big\langle\g + \zeta, \wh H_{p\rho_2}(\xi,\tilde \xi)[\tilde \g + \tilde \zeta]\Big]\Big\rangle\bigg]\\
&&\dis \le \tilde\dbE\bigg[\Big\langle\g + \zeta,  \wh H_{pp}(\xi)[\g + \zeta]\Big\rangle+\Big\langle \g + \zeta,  \Xi+2\l \wh H_{pp}(\xi)\eta\Big\rangle+{c_0\over 2}[|\g + \zeta|^2+|\tilde \g + \tilde \zeta|^2]\Big]\bigg]\\
&&\dis = \dbE\bigg[\Big\langle\g + \zeta,  [\wh H_{pp}(\xi)+ c_0I_d][\g + \zeta]\Big\rangle+\Big\langle \g + \zeta,  \tilde \dbE_{\cF^1_T}[\Xi]+2\l \wh H_{pp}(\xi)\eta\Big\rangle\bigg]\\
&&\dis =  \dbE\bigg[- \Big|\Big[(-\wh H_{pp}(\xi)- c_0I_d)^{1\over 2}[\g + \zeta]-  {1\over 2}\L(\xi,\eta)\Big|^2 + {1\over 4} |\L(\xi,\eta)|^2\Big]\bigg]\\
&&\dis \le {1\over 4} \dbE\big[|\L(\xi,\eta)|^2\big].
\eeaa
Then clearly \reff{lambdaHH} implies \reff{Hdisplacement}.
\qed

\begin{rem}
\label{rem-displacement} For standard MFGs where $b, f$ do not depend on the law of $\a$, we have $\wh H(x, p, \rho) = H(x, p, \mu)$ where $\mu = {\pi_1}_\# \rho$, and thus $\pa_{\rho_1} \wh H = \pa_\mu H$, $\pa_{\rho_2} \wh H = 0$, $c_0=0$. 
Note that $H$ is concave in $p$. We shall assume it is strictly concave and thus $H_{pp} <0$. Then \reff{lambdaHH} reduces to 
\bea
\label{Hdisplacement3}
\left.\ba{c}
\dis \tilde\dbE\bigg[ \Big\langle \eta,  [ H_{xx}(\xi)-2\l  H_{px}(\xi)]\eta\Big\rangle + \Big\langle \eta,  ~[ H_{x\mu}(\xi,\tilde \xi)-2\l H_{p\mu}(\xi,\tilde \xi)] \tilde \eta\Big\rangle \\
\dis - {1\over 4} \Big|(-H_{pp}(\xi))^{-{1\over 2}}\big[\tilde \dbE_{\cF^1_T}[H_{p\mu}(\xi, \tilde \xi)\tilde\eta]+2\l  H_{pp}(\xi)\eta\big]\Big|^2\bigg]\geq 0.
\ea\right.
\eea
Moreover, when $\l=0$, \reff{Hdisplacement3} (and \reff{lambdaHH2}) reduces further to 
\bea
\label{Hdisplacement4}
\dis \tilde\dbE\bigg[ \Big\langle \eta,  H_{xx}(\xi)\eta\Big\rangle + \Big\langle \eta,  ~H_{x\mu}(\xi,\tilde \xi) \tilde \eta\Big\rangle  - {1\over 4} \Big|(-H_{pp}(\xi))^{-{1\over 2}}\tilde \dbE_{\cF^1_T}[H_{p\mu}(\xi, \tilde \xi)\tilde\eta]\Big|^2\bigg]\geq 0.
\eea
This is exactly the condition in \cite[Definition 3.4]{GMMZ}, except that \cite{GMMZ} uses $-H$ instead of $H$. 
\end{rem}

We now present an example which satisfies \reff{lambdaHH2}, and hence \reff{Hdisplacement} with $\l=0$.

\begin{eg}
\label{eg-displacement}
We consider a special case of \reff{nonseparable} with $d=1$: for some constant $0<c<1$,
\bea
\label{nonseparable3}
\left.\ba{c}
\dis b(x, a, \cL_{(\xi,\a)}) = -a  + b_1(\cL_\xi, \dbE[\a]), \q
\dis f(x, a, \cL_{(\xi,\a)}) = {|a|^2\over 2} -  c a \dbE[\a]  + f_1(x, \cL_{(\xi, \a)}).
\ea\right.
\eea
Assume there exist constants $0\le c_0<1$ and  $\k>0$ such that
\bea
\label{lambdaHH3}
|\bar c \pa_{m_2} b_1 -\hat c|\le c_0,\q \pa_{xx} f_1 \ge \k \ge \|\pa_{x \nu_1} f_1\| + {1\over 4(1-c_0)}\Big(\|\pa_{m_1} b_1\|+[1+ \hat c]\|\pa_{x \nu_2} f_1\|\Big)^2,
\eea
where $\hat c:= {c\over 1-c}$, $\bar c:=\frac{1}{1-c}$, $m_1, m_2$ stand for $\dbE[\xi]$ and $\dbE[\a]$, respectively, and $\|\cd\|$ denotes the supremum norm of the function over all variables. Then \reff{lambdaHH2} holds true.
\end{eg}
\proof By Example \ref{eg-nonseparable} we see that
\beaa
&\dis \Phi(\cL_{(\xi, \eta)}) = \cL_{(\xi, ~ \hat c \dbE[\eta]  + \eta)},\\
&\dis H(x, p, \cL_{(\xi,\a)}) = -{1\over 2} \Big|c\dbE[\a]  + p\Big|^2 + p b_1(\cL_\xi, \dbE[\a])  +f_1(x, \cL_{(\xi,\a)}),\\
&\dis \wh H(x, p, \cL_{(\xi, \eta)}) =  -{1\over 2} \Big|\hat c\dbE[\eta]   +p\Big|^2 +p b_1(\cL_\xi, \bar c\dbE[\eta]) + f_1(x,  \cL_{(\xi, ~ \hat c \dbE[\eta]  + \eta)}).
\eeaa
Applying Proposition \ref{prop-chain} and Example \ref{eg-derivative} we have, for $\wh f_1(x,\rho):= f_1(x, \Phi(\rho))$ where $\rho=\cL_{(\xi, \eta)}$,
\beaa
&\dis \pa_{\rho_1} \wh f_1(x, \rho, \tilde x, \tilde p) = \pa_{\nu_1} f_1(x, \Phi(\rho), \tilde x, \tilde p + \hat c \dbE_{\rho}[\eta]),\\
&\dis \pa_{\rho_2} \wh f_1(x, \rho, \tilde x, \tilde p) = \pa_{\nu_2} f_1(x, \Phi(\rho), \tilde x, \tilde p + \hat c \dbE_{\rho}[\eta]) + \hat c \dbE_{\rho}\Big[ \pa_{\nu_2} f_1(x, \Phi(\rho), \xi, \eta + \hat c \dbE_{\rho}[\eta])\Big].
\eeaa
Then one may compute straightforwardly that
\beaa
&\wh H_{pp} = -1,\q  \wh H_{xx} = \pa_{xx} f_1, \q \wh H_{x\rho_1} = \pa_{x \nu_1} f_1,\q  \wh H_{p \rho_1}= \pa_{m_1} b_1,\\
&\wh H_{x\rho_2} =  \pa_{x\nu_2} f_1 + \hat c \dbE_{\rho}[ \pa_{x\nu_2} f_1],\q \wh H_{p\rho_2}= \bar c \pa_{m_2} b_1 -\hat c.
\eeaa
Then $|\wh H_{p\rho_2}|\le c_0$, and \reff{lambdaHH2} becomes
\bea
\label{lambdaHH4}
\left.\ba{c}
\dis \tilde\dbE_{\rho}\Big[ \pa_{xx} f_1|\eta|^2  + \pa_{x\nu_1}f_1 \eta\tilde\eta - {|\L(\xi)|^2\over 4(1-c_0)}\Big] \ge 0,\q \mbox{where}\\
\dis  \L(x) :=  \bar{\tilde\dbE}_\rho\bigg[\tilde \eta \Big[\pa_{m_1} b_1\big(\cL_\xi, \tilde \xi, \bar c \dbE_{\rho}[\eta]\big)  + \pa_{x\nu_2} f_1\big(\tilde \xi, \Phi(\rho), x, \f(x) + \hat c \dbE_{\rho}[\f(\xi)]\big) \\
\dis + \hat c  \pa_{x\nu_2} f_1(\tilde \xi, \Phi(\rho), \bar \xi, \f(\bar \xi) + \hat c \dbE_{\rho}[\f(\xi)])\Big]\bigg].
\ea\right.
\eea
Clearly \reff{lambdaHH3} implies \reff{lambdaHH4}, and hence \reff{lambdaHH2}.
\qed

\subsection{Global wellposedness for master equations of standard MFGs}
\label{sect-wellposedness}

For standard MFGs, by combining Proposition \ref{prop-dissufficient} and the strategy in \cite{GMMZ}, see also Subsection \ref{sect-strategy}, one can easily establish the following global wellposedness result  for the master equation under displacement semi-monotonicity, which generalizes \cite[Theorem 6.3]{GMMZ}. We remark again that, for MFGC master equations, we shall investigate their global wellposedness in future research.

\begin{thm}
\label{thm-global}

Assume $\l\geq 0$, $b(x, a, \nu)=a$ and $f(x, a, \nu) = f(x, a, \mu)$. Assume further that:

(i) $H$ and $G$ have the regularity:
\beaa
&\dis H,\pa_{xx}H,\pa_{xp}H,\pa_{pp}H,\pa_{xxp}H,\pa_{xpp}H,\pa_{ppp}H\in \cC^2(\dbR^{2d}\times\cP_2(\dbR^d)),\\
&\dis \pa_{x\mu}H,\pa_{p\mu}H,\pa_{xp\mu}H,\pa_{pp\mu}H\in \cC^2(\dbR^{2d}\times\cP_2(\dbR^d)\times\dbR^d),\\
&\dis G,\pa_{xx}G\in \cC^2(\dbR^d\times\cP_2(\dbR^d)),\qq \pa_{x\mu}G\in\cC^2\dbR^d\times\cP_2(\dbR^d)\times\dbR^d,
\eeaa
and all the second and higher order derivatives of $H$ and $G$ involved above are uniformly bounded;

(ii) $H$ is uniformly concave in $p$: $\pa_{pp}H\leq -c_0 I_d$ for some constant  $c_0>0$.

(iii)  \eqref{Hdisplacement3} holds for $H$ and \eqref{lambdadisplacementsm} holds for $G$.

\no Then the master equation \eqref{master} on $[0,T]$ admits a unique classical solution $V$ with bounded $\pa_{xx}V$ and $\pa_{x\mu}V$.
\end{thm}
\proof We shall follow the road map given in Section 3.2 to show the global wellposedness. Since the arguments are very similar to those in \cite{GMMZ, MZ3}, at below we sketch a proof only.

{\it Step 1.} We apply Theorem \ref{thm:displacement} to show that, if $V$ is a classical solutions of the master equation \eqref{master} and $V$ satisfies the Assumptions \ref{assum-regV}, then $V$ propagates the displacement $\l$-monotonicity, i.e. $V(t,\cdot,\cdot)$ satisfies \eqref{lambdadisplacementsm} for all $t\in[0,T]$.

{\it Step 2.} We shall follow the same proof as the one in \cite[Theorem 5.1]{GMMZ} to show an a priori uniform $\cW_2$-Lipschitz continuity of $\pa_xV$ in $\mu$, uniformly in $(t,x)\in[0,T]\times\dbR^d$. We note that $V$ might not be uniformly $\cW_2$-Lipschitz continuous in $\mu$ under our (weaker) assumptions. The key assumption we used in \cite[Theorem 5.1]{GMMZ} is the boundedness of $\pa_{xx}V$, which was proved using the first order derivatives of $H$ and $G$ in \cite[Proposition 6.1]{GMMZ}. This is not the case anymore here. To show it, we first apply Theorem \ref{thm:displacement} to prove that $V(t,\cd,\cd)$ satisfies \reff{lambdadisplacementsm} for all $t\in[0,T]$. By Remark \ref{rem-mon}-(iv), $\pa_{xx}V$ is uniformly semi-convex in $x$, uniformly in $(t,\mu)\in [0,T]\times\cP_2(\dbR^d)$. It is standard to obtain the uniform semi-concavity of $V$ in $x$ from the boundedness of the second order derivatives of $H$ and $G$ by the classical control theory. Thus, we obtain the a priori boundedness of $\pa_{xx}V$. Then we obtain the uniform $\cW_2$-Lipschitz continuity of $V$ in $\mu$. 

By \cite[Proposition 6.2]{GMMZ}, we can further strengthen the above a priori $\cW_2$-Lipschitz continuity to an a priori $\cW_1$-Lipschitz continuity for $\pa_{x}V$ in $\mu$.

{\it Step 3.} We shall follow the same proof as the one in \cite[Theorem 7.1]{MZ3} to show the global wellposedness of the master equation \eqref{master}. The desired regularity of the solution $V$ is a byproduct of {\it Step 2}. However, we cannot  show directly the wellposedness of the master equation due to the lack of the a priori Lipschitz continuity of $V$ in $x$ and $\mu$, we thus use the approach in \cite[Section 7]{MZ3}. That is, we first use the a priori Lipschitz estimate of $\pa_xV$ constructed in {\it Step 2} to show the wellposedness of the vectorial master equation for $\vec U := \pa_x V$. We then utilize the solution to the vectorial master equation to establish the wellposedness of the master equation \eqref{master}. 
\qed

\begin{rem}
If $G$ satisfies the Lasry-Lions monotonicity and $\pa_{xx}G$ is bounded by $\l$, then $G$ is displacement semi-monotone. Therefore,  we obtain that, if $H$ and $G$ satisfy the assumptions (i), (ii) in Theorem \ref{thm-global}, $H$ satisfies \eqref{Hdisplacement3}  and $G$ is Lasry-Lions monotone, then the master equation is wellposed on $[0,T]$. In this sense, Theorem \ref{thm-global} unifies the wellposedness results under the Lasry-Lions monotonicity and the displacement monotonicity.

We shall remark though, even when $G$ is Lasry-Lions monotone,  $V$ propagates the displacement semi-monotonicity, not necessarily the Lasry-Lions monotonicity (when $f$ is non-separable).
\end{rem}

\section{Propagation of anti-monotonicity}
\label{sect-anti}
\setcounter{equation}{0}
In this section we fix $\vec \l \in D_4$. Recall \reff{kappaA}.

\begin{assum}\label{assum-antidisplacement} 
(i) $\wh H\in \cC^2(\mathbb R^{d}\times \dbR^d\times\mathcal{P}_2(\mathbb R^{2d}))$  and there exist constants $\ol L, L_0>0$ and  $\overline \g >\underline \g  >0$  such that 
\bea
\label{Hbound}
&\dis |\pa_{xp}\wh H|\leq \overline \g L_0,\quad |\pa_{xx}\wh H|\leq \overline\g L_0,\q
|\pa_{pp}\wh H|,|\pa_{x\rho_1}\wh H|, |\pa_{x\rho_2}\wh H|, |\pa_{p\rho_1}\wh H|, |\pa_{p\rho_2}\wh H|\le \ol L;\\
\label{antiH1}
&\dis \underline\k(-\pa_{xp}\wh H)\geq L_0,\q  \underline\k(-\pa_{xx}\wh H) \geq\underline \g L_0.
\eea

(ii) There exists a constant $L^u_{xx}>0$ such that 
\bea
\label{theta}
 \th_1 := { \overline \g [1+L^u_{xx}]  \over \sqrt{4( \underline \g \l_0+2\l_3)}} <1,\q\mbox{and}\q \ol L \ol \k(A_1^{-1}A_2) \le L_0,
\eea
where:
\bea
\label{A}
&&\dis A_1:= 
\begin{bmatrix}
4[1-\th_1]   & 0& 0 \\
0 &  2\l_2 & 0\\
0&  0& [1-\th_1][\l_0\underline \g+2\l_3 ] 
\end{bmatrix},\nonumber\\
&&\dis A_2:= B_1L_{xx}^u+B_2:= \begin{bmatrix}
2 & 2+\l_2 & 1\\
2+\l_2 &  4\l_2 &\l_2\\
1 &  \l_2& 0
\end{bmatrix} L_{xx}^u   \\
&&\dis +\begin{bmatrix}
 \l_0+2|\l_0-\l_1| & \l_0+|\l_0-\frac{1}{2}\l_1|+\frac{1}{2}|\l_1|+\l_2 & |\l_0-{1\over 2}\l_1|+\frac{1}{2}|\l_1|+2\l_3 \\
\l_0+|\l_0-\frac{1}{2}\l_1|+\frac{1}{2}|\l_1|+\l_2  &  2| \l_1|+2\l_2 &|\l_1|+ \l_2 +2\l_3 \\
|\l_0-{1\over 2}\l_1|+\frac{1}{2}|\l_1|+2\l_3  & |\l_1|+ \l_2 +2\l_3  &  |\l_1|+2\l_3
\end{bmatrix}.\nonumber
\eea
\end{assum}

\ms

\begin{thm}\label{thm:anti}
Let Assumptions \ref{assum-fix},  \ref{assum-regH}, \ref{assum-regV}, and  \ref{assum-antidisplacement} hold. Assume further that, for the constant  $L^u_{xx} $ in Assumption \ref{assum-antidisplacement} (ii),
\bea
\label{Vxx}
|\pa_{xx} V|\le L^u_{xx}.
\eea
 If $G$ satisfies the $\vec \l$-anti-monotonicity \reff{anti}, then $V(t,\cd,\cd)$ satisfies \reff{anti} for all $t\in[0,T]$. 
\end{thm}

We remark that the bound $L^u_{xx}$ of $\pa_{xx} V$ can be estimated a priori by using the HJB equation or the backward SDE in the mean field game system, see \cite[Section 6]{MZ3} for more details.

\proof 
Without loss of generality, we shall prove the theorem only for $t=0$. We will continue to use the notation in the proofs of Theorem \ref{thm:LL} and \ref{thm:displacement}. Introduce:
\beaa
\Xi_t:=\l_0\bar I(t)+\l_1 I(t)+\mathbb E\big[|\Gamma_t|^2+\l_2|\Upsilon_t|^2-\l_3|\delta X_t|^2\big] = MON^{anti}_{\vec \l} V(t,\cd, \cd)(X_t, \d X_t).
\eeaa
Then it is sufficient to show that
\bea
\label{dXi}
{d\over dt} \Xi_t \ge 0.
\eea

Following the calculation in \cite[Theorem 4.1]{MZ3}  we have
\bea\label{dUpsilon}
\left.\ba{lll}
\dis d\Upsilon_t = \big[ -K_1(t) \Upsilon_t - K_2(t)\big] dt  + (dB_t)^\top K_3(t) + \b  (dB^0_t)^\top K_4(t);\\
d\G_t = \big[ -2\wh H_{xp}(X_t) \G_t+\pa_{xx}V(X_t)\wh H_{pp}(X_t)\Upsilon_t - \bar K_1(t)\big] dt + (dB_t)^\top \bar K_2(t) + \b  (dB^0_t)^\top \bar K_3(t),
\ea\right.
\eea
where 
\bea
\label{K}
\left.\ba{lll}
K_1(t) := \wh H_{xp}(X_t) + \pa_{xx} V(X_t) \wh H_{pp}(X_t),\\
  K_2(t) := \tilde\dbE_{\mathcal{F}_t}\Big[  \big[\wh H_{x\rho_1}(X_t,\tilde X_t) \delta\tilde X_t+\wh H_{x\rho_2}(X_t,\tilde X_t)[\tilde \Gamma_t+\tilde \Upsilon_t]\big] \\
  \qq\qq\qq\q+  \pa_{xx} V(X_t)\big[\wh H_{p\rho_1}(X_t,\tilde X_t) \delta\tilde X_t+\wh H_{p\rho_2}(X_t,\tilde X_t)[\tilde \Gamma_t+\tilde \Upsilon_t]\big]\Big], \\
K_3(t) := \tilde{\mathbb E}_{\mathcal{F}_t}\big[\pa_{xx\mu}V(X_t,\tilde X_t)\delta \tilde X_t\big],\\
 K_4(t):= K_3(t) + \bar{\tilde \dbE}_{\mathcal{F}_t}\Big[\big[(\pa_{\mu x\mu}V)(X_t,\bar X_t,\tilde X_t)  +\pa_{\tilde xx\mu}V(X_t,\tilde X_t)\big]\delta\tilde X_t\Big],\\
 \bar K_1(t) :=[\wh H_{xx}(X_t)-\pa_{xx}V(X_t)\wh H_{px}(X_t)]\delta X_t\\
 \qq\qq- \pa_{xx}V(X_t)\tilde\dbE_{\mathcal{F}_t} \big[\wh H_{p\rho_1}(X_t,\tilde X_t)\delta \tilde X_t+\wh H_{p\rho_2}(X_t,\tilde X_t)[\tilde \Gamma_t+\tilde \Upsilon_t]\big], \\
 \bar K_2(t) := \pa_{xxx}V(X_t)\delta  X_t,\\
 \bar K_3(t):= \bar K_2(t) +  \tilde \dbE_{\mathcal{F}_t}\Big[(\pa_{\mu xx}V)(X_t,\tilde X_t)\delta\tilde X_t\Big].
\ea\right.
\eea
In particular, this implies that
\bea\label{dUpsilon2}
\left.\ba{lll}
\dis {d\over dt} \dbE[|\Upsilon_t|^2]  \ge  2\dbE\Big[\big\langle \Upsilon_t,~  -K_1(t) \Upsilon_t - K_2(t)\big\rangle\Big];\ss\\
\dis {d\over dt} \dbE[|\G_t|^2] \ge 2 \dbE\Big[\big\langle \G_t ,~ -2\wh H_{xp}(X_t) \G_t+\pa_{xx}V(X_t)\wh H_{pp}(X_t)\Upsilon_t - \bar K_1(t)\big\rangle\Big].
\ea\right.
\eea
Thus, combining \reff{dtI-LL}, \reff{cdbarI},  and \reff{deltaX2}, and recalling the $N$ in \reff{XY}, we have
\beaa
\left.\ba{lll}
\dis  {d\over dt}\Xi_t  \ge \l_0 \tilde\dbE\Big[\big\langle \wh H_{pp}(X_t)\G_t,\G_t\big\rangle+2\big\langle \wh H_{pp}(X_t)\G_t, \Upsilon_t\big\rangle\\
\dis\qq\qq + 2\big\langle\G_t, \wh H_{p\rho_1}(X_t,\tilde X_t)\delta \tilde X_t+\wh H_{p\rho_2}(X_t,\tilde X_t)[\tilde\Gamma_t+\tilde \Upsilon_t]\big\rangle-\big\langle \wh H_{xx}(X_t)\delta X_t,\delta X_t\big\rangle \Big] \\
\dis ~ + \l_1  \tilde\dbE\Big[\big\langle \wh H_{pp}(X_t) \Upsilon_t,~\Upsilon_t\big\rangle - \big\langle \wh H_{x\rho_1}(X_t,\tilde X_t) \delta  \tilde X_t+\wh H_{x\rho_2}(X_t,\tilde X_t)[\tilde\G_t+\tilde\Upsilon_t], \delta  X_t\big\rangle\\
\dis\qq\qq -\big\langle  \wh H_{p\rho_1}(X_t,\tilde X_t)\delta  \tilde X_t+\wh H_{p\rho_2}(X_t,\tilde X_t)[\tilde\G_t+\tilde \Upsilon_t], \G_t-\Upsilon_t \big\rangle  \Big]\\
\dis ~+2\tilde\dbE\Big[ \big\langle \G_t, \big[ -2\wh H_{xp}(X_t) \G_t+\pa_{xx}V(X_t)\wh H_{pp}(X_t)\Upsilon_t - \bar K_1(t)\big] \big\rangle +\l_2\big\langle\Upsilon_t, \big[ -K_1(t) \Upsilon_t - K_2(t)\big] \big\rangle\Big]\\
\dis~ -2\l_3\tilde\dbE\Big[\big\langle \wh H_{px}(X_t)\delta X_t+\wh H_{p\rho_1}(X_t,\tilde{X_t})\delta \tilde X_t+\wh H_{p\rho_2}(X_t,\tilde X_t)[\tilde \G_t+\tilde\Upsilon_t]+\wh H_{pp}(X_t)[ \Upsilon_t+\G_t],\delta X_t \big\rangle\Big]\\
\dis = \tilde\dbE\bigg[\Big\langle \big[ \l_0 \wh H_{pp}(X_t) -4 \wh H_{xp}(X_t)\big]  \G_t, ~\G_t\Big\rangle + \Big\langle \big[\l_1\wh H_{pp}(X_t)- 2\l_2K_1(t)] \Upsilon_t, \Upsilon_t \Big\rangle\\
\dis ~ + \Big\langle \big[-\l_0 \wh H_{xx}(X_t) -2\l_3\wh H_{px}(X_t)\big] \d X_t, \d X_t\Big\rangle +\Big\langle \big[2\l_0-\l_1+2\pa_{xx}V(X_t)\big]\wh H_{p\rho_2}(X_t,\tilde X_t)\tilde \Gamma_t, \G_t \Big\rangle\\
\dis ~ +\Big \langle\big[[\l_1-2\l_2\pa_{xx}V(X_t)]\wh H_{p\rho_2}(X_t,\tilde X_t)-2\l_2\wh H_{x\rho_2}(X_t,\tilde X_t)\big]\tilde \Upsilon_t, \Upsilon_t\Big\rangle \\
\dis~ -\Big\langle \big[ \l_1\wh H_{x\rho_1}(X_t,\tilde X_t) +2\l_3\wh H_{p\rho_1}(X_t,\tilde X_t)\big]\delta \tilde X_t, \d X_t\Big\rangle\\
\dis~ +\Big\langle  2 [\l_0 \wh H_{pp}(X_t) + \pa_{xx} V(X_t) \wh H_{pp}(X_t)]\Upsilon_t+[2\l_0-\l_1+2\pa_{xx}V(X_t)]\wh H_{p\rho_2}(X_t,\tilde X_t)\tilde \Upsilon_t,\G_t\Big\rangle \\
\dis~+\Big\langle \big[[ \l_1-2\l_2\pa_{xx}V(X_t)]\wh H_{p\rho_2}(X_t,\tilde X_t) -2\l_2\wh H_{x\rho_2}(X_t,\tilde X_t)\big]\tilde \Gamma_t,\Upsilon_t\Big\rangle\\
\dis~   +\Big\langle \big[ [2\l_0-\l_1+2\pa_{xx}V(X_t)]  \wh H_{p\rho_1}(X_t,\tilde X_t) -\l_1\wh H_{\rho_2x}(\tilde X_t,X_t)-2\l_3\wh H_{\rho_2p}(\tilde X_t,X_t)\big]\delta\tilde X_t\\
\dis~\q\,\,+2\big[-\wh H_{xx}(X_t)+\pa_{xx}V(X_t)\wh H_{px}(X_t)-\l_3\wh H_{pp}(X_t)\big]\delta X_t, \G_t\Big\rangle \\
\dis~ + \Big\langle  \big[[\l_1-2\l_2\pa_{xx}V(X_t)]  \wh H_{p\rho_1}(X_t,\tilde X_t) -\l_1\wh H_{\rho_2x}(\tilde X_t,X_t)\\
 \dis~\q-2\l_2\wh H_{x\rho_1}(X_t,\tilde X_t)-2\l_3\wh H_{\rho_2 p}(\tilde X_t, X_t)\big]\delta \tilde X_t-2\l_3\wh H_{pp}(X_t)\delta X_t ,\Upsilon_t\Big\rangle \bigg].
\ea\right.
\eeaa
Recall \reff{D4}, \reff{kappaA}, and \reff{K},   by \reff{Hbound} and \reff{antiH1} we have
\beaa
\left.\ba{lll}
\dis  {d\over dt}\Xi_t  \ge \big[ 4 L_0 - \l_0\ol L\big] \dbE[|\G_t|^2]  + \big[ 2\l_2 L_0 - |\l_1|\ol L - 2\l_2 L^u_{xx} \ol L\big] \dbE[|\Upsilon_t|^2] \\
\dis \qq +  \big[\l_0 \ul \g L_0 + 2\l_3  L_0\big] \dbE[|\d X_t|^2]  -\big[ |2\l_0-\l_1| + 2L^u_{xx}\big] \ol L \big(\dbE[|\G_t|]\big)^2 \\
\dis\qq -\big[|\l_1|+2\l_2L^u_{xx} +2\l_2\big]\ol L  \big(\dbE[|\Upsilon_t|]\big)^2 -\big[ |\l_1|  +2\l_3\big]\ol L\big(\dbE[|\d X_t|]\big)^2\\
\dis\qq -\big[2 [\l_0 + L^u_{xx}]  +|2\l_0-\l_1| + 2L^u_{xx} + |\l_1|+2\l_2 L^u_{xx} +2\l_2 \big]\ol L \dbE[|\G_t|] \dbE[| \Upsilon_t|] \\
\dis\qq  - \big[ [|2\l_0 -\l_1|+2L^u_{xx}  + |\l_1| + 2\l_3]  \ol L  +2[\ol \g L_0  + L^u_{xx} \ol \g L_0 + \l_3 \ol L]\big] \dbE[|\d X_t|] \dbE[|\G_t|]  \\
\dis\qq - \big[[|\l_1|+2\l_2L^u_{xx} + |\l_1| + 2\l_2 +2\l_3 +2\l_3] \ol L  \dbE[|\delta X_t|]\dbE[|\Upsilon_t|]\\
\dis\q \ge \Big[ 4 L_0 - \big[\l_0+ |2\l_0-\l_1| + 2L^u_{xx}\big] \ol L\Big] \big(\dbE[|\G_t|]\big)^2  \\
\dis \qq + \Big[ 2\l_2 L_0 - \big[ 2|\l_1| + 4\l_2 L^u_{xx}+2\l_2\big] \ol L\Big]   \big(\dbE[|\Upsilon_t|]\big)^2 \\
\dis \qq +  \Big[[\l_0 \ul \g +  2\l_3 ]L_0  -\big[ |\l_1|  +2\l_3\big]\ol L\Big]\big(\dbE[|\d X_t|]\big)^2\\
\dis\qq - 2\Big[\l_0 +|\l_0-{\l_1\over 2}| + {|\l_1|\over 2}+ \l_2 + [2+\l_2] L^u_{xx}  \Big]\ol L \dbE[|\G_t|] \dbE[| \Upsilon_t|] \\
\dis\qq  - 2\Big[ \overline \g  [1+ L^u_{xx}] L_0  + \big[|\l_0-{\l_1\over 2}| + {|\l_1|\over 2}+2\l_3+L^{u}_{xx}]\ol L\Big] \dbE[|\d X_t|] \dbE[|\G_t|]  \\
\dis\qq - 2\Big[|\l_1|+ \l_2 +\l_2 L^u_{xx} +2\l_3\Big] \ol L  \dbE[|\delta X_t|]\dbE[|\Upsilon_t|]\\
\dis\q \ge \Big[ 4 [1-\th_1]L_0 - \big[\l_0+ |2\l_0-\l_1| + 2L^u_{xx}\big] \ol L\Big] \big(\dbE[|\G_t|]\big)^2  \\
\dis \qq + \Big[ 2\l_2 L_0 - \big[ 2|\l_1| + 4\l_2 L^u_{xx}+2\l_2\big] \ol L\Big]   \big(\dbE[|\Upsilon_t|]\big)^2 \\
\dis \qq +  \Big[(1-\th_1)[\l_0 \ul \g +  2\l_3 ]L_0  -\big[ |\l_1|  +2\l_3\big]\ol L\Big]\big(\dbE[|\d X_t|]\big)^2\\
\dis\qq - 2\Big[\l_0 +|\l_0-{\l_1\over 2}| + {|\l_1|\over 2}+ \l_2 + [2+\l_2] L^u_{xx}  \Big]\ol L \dbE[|\G_t|] \dbE[| \Upsilon_t|] \\
\dis\qq  - 2\Big[  \big[|\l_0-{\l_1\over 2}| + {|\l_1|\over 2}+2\l_3+L^{u}_{xx}]\ol L\Big] \dbE[|\d X_t|] \dbE[|\G_t|]  \\
\dis\qq - 2\Big[|\l_1|+ \l_2 +\l_2 L^u_{xx} +2\l_3\Big] \ol L  \dbE[|\delta X_t|]\dbE[|\Upsilon_t|],
\ea\right.
\eeaa
where in the last step we used the fact that: recalling the $\th_1$ in \reff{theta},
\beaa
  2\overline \g  [1+ L^u_{xx}] \dbE[|\d X_t|] \dbE[|\G_t|]  \le 4 \th_1 \big(\dbE[|\G_t|]\big)^2  + \th_1 [\l_0\underline \g +2\l_3] \big(\dbE[|\d X_t|]\big)^2,
\eeaa
Then, recalling \reff{A} and denoting $e := \big( \dbE[|\G_t|],~ \dbE[|\Upsilon_t|],~   \dbE[|\d X_t|]\big)$, we have
\beaa
 {d\over dt} \Xi_t \ge    e \big[ A_1  L_0 - A_2 \ol L\big] e^\top \ge 0,
\eeaa
thanks to  \reff{theta} and the fact that $A_1 > 0$.
\qed

\begin{eg}
\label{eg-anti}
Again we consider a special case of \reff{nonseparable} with $d=1$: 
\beaa
\left.\ba{c}
\dis b(x, a, \cL_{(\xi,\a)}) = -a  -L_0 x+ b_1(\cL_\xi, \dbE[\a]), ~
\dis f(x, a, \cL_{(\xi,\a)}) = {|a|^2\over 2} -  c a \dbE[\a]  - {\gamma L_0\over 2} x^2+ f_1(x, \cL_{(\xi, \a)}),
\ea\right.
\eeaa
for some constants $0<c<1$, $\gamma>0$, $L_0>0$. For any $L^u_{xx}>0$, when $L_0$ is large enough, there exist appropriate $\ol L>0$, $\ol \g > \ul \g >0$, and $\vec{\l}\in D_4$ such that Assumption \ref{assum-antidisplacement} holds true.
\end{eg}
\proof By Example \ref{eg-nonseparable} and recalling the notations $\hat c, \bar c$ in Example \ref{eg-displacement}, we see that 
\beaa
&\dis \Phi(\cL_{(\xi, \eta)}) = \cL_{(\xi, ~ \hat c \dbE[\eta]  + \eta)},\\
&\dis H(x, p, \cL_{(\xi,\a)}) = -{1\over 2} \Big|c\dbE[\a]  + p\Big|^2 - L_0 xp+ p b_1(\cL_\xi, \dbE[\a]) - {\gamma L_0\over 2} x^2 +f_1(x, \cL_{(\xi,\a)}),\\
&\dis \wh H(x, p, \cL_{(\xi, \eta)}) =  -{1\over 2} \Big|\hat c\dbE[\eta]   +p\Big|^2 -L_0 xp +p b_1(\cL_\xi,\bar c\dbE[\eta])- {\gamma L_0\over 2} x^2 + f_1(x,  \cL_{(\xi, ~ \hat c \dbE[\eta]  + \eta)}).
\eeaa
Following the same calculation as that in Example \ref{eg-displacement},
\beaa
&\wh H_{pp} = -1,\q \wh H_{xp}=- L_0,\q  \wh H_{xx} =  -\gamma L_0 +\pa_{xx}f_1, \q \wh H_{x\rho_1} = \pa_{x \nu_1} f_1,\\
&\wh H_{x\rho_2} =  \pa_{x\nu_2} f_1 + \hat c \dbE[ \pa_{x\nu_2} f_1],\q  \wh H_{p \rho_1}= \pa_{m_1} b_1,\q \wh H_{p\rho_2}= \bar c \pa_{m_2} b_1 -\hat c.
\eeaa

For given functions $f_1, b_1$, clearly there exists a fixed constant $\ol L>0$ such that 
\[
|\pa_{pp}\wh H|,|\pa_{x\rho_1}\wh H|, |\pa_{x\rho_2}\wh H|, |\pa_{p\rho_1}\wh H|, |\pa_{p\rho_2}\wh H|\le \ol L. 
\]
Set $\ul \g := {\g\over 2}$, $\ol \g := \g + 1$, for $L_0$ sufficiently large, we have
\beaa
&\ul \k (-\wh H_{xp}) = L_0,\qq  |\wh H_{xp})| = L_0 \le \ol \k L_0;\\
&|\pa_{xx} \wh H| - \ol \g L_0 \le \g L_0 + \ol L - \ol \g L_0 = \ol L - L_0 \le 0;\\
& \ul\k (-\pa_{xx} \wh H) - \ul \g L_0 \ge \g L_0 - \ol L - \ul \g L_0 = {\g\over 2} L_0 - \ol L \ge 0.
\eeaa
That is,  \eqref{Hbound} and \eqref{antiH1} hold true.

We next fix arbitrary $\l_0, \l_2 >0$ and $\l_1 \in \dbR$. Choose $\l_3> 0$ sufficiently large we have $\th_1 \le {1\over 2}$. Finally, set $L_0$ sufficiently large such that $L_0 \ge \ol L \ol \k(A_1^{-1}A_2)$, we verify \reff{theta} as well.
\qed

We shall point out though, by \reff{Vxx} $L^u_{xx}$ may in turn depend on $L_0$, so extra efforts are needed in order to ensure full compatibility of our conditions. This, however, requires the a priori estimate for $\pa_{xx} V$ which is not carried out in this paper. We thus leave it to our accompanying paper on global wellposedness of MFGC master equations. We remark that we have a complete result in \cite{MZ3} for standard MFG master equations. 

\section{Appendix}
\label{sect-proof}
\setcounter{equation}{0}

\no {\bf Proof of \reff{dtI-LL}}. We first apply the It\^o's formula \reff{Ito} on $\pa_{x\mu}V(t, X_t, \cL_{X_t|\cF^0_t}, \tilde X_t)$ to obtain 
\begin{equation}\label{eq:ito_double1}
 {d\over dt} I(t)= I_1 + I_2 + I_3,
\end{equation}
where, by using $\hat X$ to denote another conditionally independent copy, 
\beaa
I_1&:=& \bar{\tilde \dbE}\bigg[\bigg\langle \Big\{\pa_{tx\mu} V(X_t, \tilde X_t)  + {\widehat \beta^2\over 2} ((\tr\pa_{xx})\pa_{x\mu} V)(X_t,  \tilde X_t)  +\wh H_p (X_t)^\top\pa_{xx\mu} V(X_t,  \tilde X_t)\\
&&+\beta^2(\tr(\pa_{x\mu})\pa_{x\mu}V)(X_t,\bar X_t,\tilde X_t)+\beta^2(\tr(\pa_{\tilde x\mu})\pa_{x\mu}V)(X_t,\bar X_t,\tilde X_t)\\
&&+\beta^2(\tr(\pa_{\tilde x x})\pa_{x\mu}V)(X_t,\tilde X_t)+\frac{\beta^2}{2}(\tr(\pa_{\mu\mu})\pa_{x\mu}V)(X_t,\hat X_t,\bar X_t,\tilde X_t) \\
&&+  {\widehat\beta^2\over 2}(\tr(\pa_{\bar x\mu}) \pa_{ x\mu} V)(X_t, \bar X_t, \tilde X_t)  + \wh H_p(\bar X_t)^\top\pa_{\mu x\mu} V(X_t, \bar X_t, \tilde X_t)\\
&& + {\widehat\beta^2\over 2} (\tr(\pa_{\tilde x\tilde x})\pa_{x\mu} V)(X_t, \tilde X_t) + \wh H_p(\tilde X_t)^\top\pa_{\tilde xx\mu} V(X_t, \tilde X_t)   \Big\}\delta  \tilde X_t, ~ \delta  X_t\bigg\rangle \bigg];\\
I_2&:=&- \hat{\bar{\tilde \dbE}}\bigg[\Big\langle \pa_{\mu x} V(X_t, \tilde X_t) \Big\{ \big[\wh H_{px} (X_t) + \wh H_{pp} (X_t)\pa_{xx} V(X_t) \big]\delta  X_t \\
&& + \Big[ \wh H_{p\rho_1}(X_t, \bar X_t)+\wh H_{p\rho_2}(X_t,\bar X_t)\pa_{xx}V(\bar X_t) + \wh H_{pp}(X_t)\pa_{x\mu} V(X_t, \bar X_t) \Big]\delta \bar X_t\\
&& + \wh H_{p\rho_2}(X_t,\bar X_t)\pa_{x\mu}V(\bar X_t,\hat X_t)\delta \hat X_t
 \Big\}, ~ \delta \tilde X_t\Big\rangle \bigg];\\
I_3&:=&- \hat{\bar{\tilde \dbE}}\bigg[\Big\langle \pa_{x\mu} V(X_t,  \tilde X_t) \Big\{ \big[\wh H_{px}(\tilde X_t)+ \wh H_{pp}(\tilde X_t)\pa_{xx} V(\tilde X_t)\big] \delta \tilde X_t \\
&& + \Big[\wh H_{p\rho_1}(\tilde X_t, \bar X_t)+\wh H_{p\rho_2}(\tilde X_t, \bar X_t)\pa_{xx}V(\bar X_t)+\wh H_{pp}(\tilde X_t) \pa_{x\mu} V(\tilde X_t,\bar X_t) \Big]\delta \bar X_t\\
&& + \wh H_{p\rho_2}(\tilde X_t,\bar X_t)\pa_{x\mu}V(\bar X_t,\hat X_t)\delta \hat X_t\Big\}, \delta  X_t\Big\rangle\bigg].
\eeaa
On the other hand, applying $\partial_{x\mu}$ to \reff{master} we obtain 
\begin{equation}\label{paxmucLV}
 0 = (\pa_{x\mu} \sL V)(t, x, \mu, \tilde x) = J_1+ J_2 +J_3,
\end{equation} 
where
\beaa
J_1&:=& \pa_{tx\mu } V (x, \tilde x)+ {\widehat \beta^2\over 2}(\tr({\pa_{xx}})\pa_{x\mu} V)(x, \tilde x)  + \wh H_{x\rho_1}(x,\tilde x) \\
&& + \wh H_{x\rho_2}(x,\tilde x)\pa_{xx}V(\tilde x)+\bar\dbE[\wh H_{x\rho_2}(x,\bar\xi)\pa_{x\mu}V(\bar \xi,\tilde x)]\\
&&+ \pa_{xx} V(x) \Big[\wh H_{p\rho_1}(x,\tilde x)+\wh H_{p\rho_2}(x,\tilde x)\pa_{xx}V(\tilde x)+ \bar\dbE[\wh H_{p\rho_2}(x,\bar\xi)\pa_{x\mu}V(\bar \xi,\tilde x)]\Big]\\
&& +\Big[\wh H_{xp}(x) +\pa_{xx} V(x) \wh H_{pp}(x) \Big] \pa_{x\mu} V(x, \tilde x) + \wh H_p(x)^\top \pa_{xx\mu} V(x, \tilde x);\\
J_2 &:=& 
{{\widehat\beta^2\over 2} (\pa_{x \tilde x } \tr(\pa_{\tilde x\mu})V)(x, \tilde x)} + \pa_{x\mu} V(x,\tilde x) \Big[\wh H_{px} (\tilde x) + \wh H_{pp}(\tilde x) \pa_{xx} V(\tilde x)\Big] \\
&&+ \wh H_p(\tilde x)^\top\pa_{\tilde x x \mu} V(x,\tilde x)+ \beta^2(\pa_{x\tilde x}\tr(\pa_{x\mu})V)(x,\tilde x)+\beta^2\bar\dbE\big[(\pa_{x\tilde x}\tr(\pa_{\mu\mu})V)(x,\bar\xi,\tilde x)\big];\\
J_3 &:=& \hat{\bar \dbE}\bigg[{\widehat \beta^2\over 2}(\tr(\pa_{\bar x\mu}) \pa_{x \mu  } V)(x, \tilde x, \bar \xi) + \wh H_p(\bar \xi)^\top \pa_{\mu x \mu} V(x,  \tilde x, \bar \xi)\\ 
&&+ \pa_{x\mu} V(x, \bar \xi) \Big[ \wh H_{p\rho_1}(\bar\xi, \tilde x)+\wh H_{p\rho_2}(\bar \xi,\tilde x)\pa_{xx}V(\tilde x)+\wh H_{p\rho_2}(\bar\xi,\hat \xi)\pa_{x\mu}V(\hat \xi,\tilde x) + \wh H_{pp} (\bar \xi) \pa_{x\mu} V(\bar \xi, \tilde x)\Big]\\
&&+\beta^2(\tr(\pa_{x\mu})\pa_{x\mu}V)(x,\tilde x,\bar \xi)+\frac{\beta^2}{2}(\tr(\pa_{\mu\mu})\pa_{x\mu}V)(x,\tilde x,\hat\xi,\bar\xi)\bigg].
\eeaa
Evaluate \eqref{paxmucLV} along $(X_t,\mu_t,\tilde X_t)$ and plug into \eqref{eq:ito_double1}. As in \cite[Theorem 4.1]{GMMZ}, by straightforward calculation we obtain \reff{dtI-LL}.
\qed

\no {\bf Proof of Remark \ref{rem-LL}}. Given $\xi_i\in \dbL^2(\cF_0)$, $i=1,2$, let $X^i$ solve the  McKean-Vlasov SDE:
\bea\label{Xi}
\left.\ba{c}
\dis X_t^i = \xi_i +\int_0^t \pa_p\wh H(X_s^i, \pa_x V(s, X_s^i, \mu_s^i), \rho_s^i) ds +  B_t+\beta B_t^0,\\
\dis \mbox{where}\q \mu_t^i:=\mathcal{L}_{X_t^i|\mathcal{F}_t^{0}},\q \rho_t^i := \cL_{(X_t^i,\pa_xV(t,X_t^i,\mu_t^i))|\mathcal{F}_t^{0}}.
\ea\right.
\eea
It is standard that the optimal control is $
\alpha^{i}_s:=\pa_p\wh H(X_s^i, \pa_x V(s, X_s^i, \mu_s^i), \rho_s^i), 
$
and thus
\bea
\label{DPP1}
 \dbE[V(t,X_t^i,\mu^i_t)] &=& \mathbb E\Big[V(t_\d ,X_{t_\d}^i,\mu^i_{t_\d})+\int_t^{t_\d} f(X_s^i,\alpha_s^{i},\rho^i_s)ds\Big],
 \eea
where $ t_\delta:=t+\delta$.
Let $\alpha^\delta$ be any admissible control in $\mathcal{A}_{t_\delta}$. Consider
\beaa
\left.\ba{c}
\dis \alpha^{i,\delta}(s,x):=\left\{\ba{lll}\pa_p\wh H(x,\pa_xV(s,x,\mu_s^i),\rho_s^i) ,\quad s\in[t,t_\delta);\\
\alpha^\delta(s,x),\qq\qq\qq\qq\, s\in [t_\delta,T].
\ea\right.\ms\\
\dis X_s^{i,\delta} = X_t^i +\int_t^s \alpha^{i,\delta}(s,X_s^{i,\delta})ds +  B_s^t+\beta B_s^{0,t},\quad s\in[t,T].
\ea\right.
\eeaa
Since $b(\cdot,a,\cdot)=a$, we have $X_s^{i,\delta}=X_s^i$ and $\alpha_s^i=\alpha^{i,\delta}(s,X_s^{i,\delta})$ for any $s\in[t,t_\delta]$. Moreover
\beaa
\dis X_s^{i,\delta} = X_{t_\delta}^i +\int_{t_\delta}^s \alpha^{\delta}(s,X_s^{i,\delta})ds +  B_s^{t_\delta}+\beta B_s^{0,t_\delta},\quad s\in[t_\delta,T].
\eeaa
Thus, for $i, j=1,2$ with $i\neq j$,
\beaa
\dis \dbE[V(t,X_t^j,\mu^i_t)]  &\leq& \mathbb E\Big[G(X_T^{j,\delta},\mu_T^i)+\int_t^Tf(X_s^{j,\delta},\alpha^{j,\delta}(s,X_s^{j,\delta}),\rho_s^i)ds\Big]\\
&=& \mathbb E\Big[G(X_T^{j,\delta},\mu_T^i)+\int_{t_\delta}^{T}f(X_s^{j,\delta},\alpha^{\delta}(s,X_s^{j,\delta}),\rho_s^i)ds+\int_t^{t_\d}f(X_s^j,\alpha_s^j,\rho^i_s)ds\Big].
\eeaa
Taking infimum over all admissible controls $\alpha^\delta$ in $\mathcal{A}_{t_\delta}$ above, we have
\bea
\label{DPP2}
\dbE[V(t,X_t^j,\mu^i_t)] &\leq& \mathbb E\Big[V(t_\d,X_{t_\d}^j,\mu^i_{t_\d})+\int_t^{t_\d}f(X_s^j,\alpha_s^j,\rho^i_s)ds\Big].
\eea
Therefore, by \reff{DPP1}, \reff{DPP2}, and \reff{fmon},
\beaa
&&\mathbb E\Big[V(t_\d,X_{t_\d}^1,\mu_{t_\d}^1)+V(t_\d,X_{t_\d}^2,\mu^2_{t_\d})-V(t_\d,X_{t_\d}^1,\mu^2_{t_\d})-V(t_\d,X_{t_\d}^2,\mu^1_{t_\d})\Big]\\
&&\q -\mathbb E\Big[V(t,X_t^1,\mu^1_t)+V(t,X_t^2,\mu^2_t)-V(t,X_t^1,\mu^2_t)-V(t,X_t^2,\mu^1_t)\Big]\\
&&\leq-\mathbb E\Big[\int_t^{t+\d} \big[f(X_s^1,\alpha_s^1,\rho^1_s)+f(X_s^2,\alpha_s^2,\rho^2_s) -f(X_s^1,\alpha_s^1,\rho^2_s)-f(X_s^2,\alpha_s^2,\rho^1_s)\big]ds\Big]\le 0.
\eeaa
Divide both sides by $\d$ and then send $\d\to 0$, we obtain
\begin{eqnarray*}
\frac{d}{dt}\dbE \Big[V(t,X_t^1,\mu^1_t)+V(t,X_t^2,\mu^2_t)-V(t,X_t^1,\mu^2_t)-V(t,X_t^2,\mu^1_t)\Big]\leq 0,
\end{eqnarray*}
which implies that, denoting $\D X_t := X^2_t - X^1_t$, 
\bea
\label{Xiest}
\frac{d}{dt}\tilde \dbE\Big[\int_{0}^1\big\langle \pa_{x\mu}V(t,X_t^1+\theta \Delta X_t,\mathcal{L}_{(X_t^1+\theta \Delta X_t)|\cF^0_t},\tilde X_t^1+\theta\Delta \tilde X_t)\Delta \tilde X_t,~\Delta X_t\big\rangle \Big]d\theta\leq 0.
\eea
Now fix $\xi, \eta\in \dbL^2(\cF_0)$ and set $\xi_1:= \xi$, $\xi_2 := \xi + \e \eta$. Then $X^1$ identifies the $X$ in \reff{XY}, and by denoting $X^\e = X^2$, one can verify that $\lim_{\e\to 0} {1\over \e}[X^\e_t - X_t] = \d X_t$, where the limit is in $\dbL^2$ sense and $\d X$ is defined in \reff{XY}. Then, by dividing \reff{Xiest} with $\e^2$ and sending $\e\to 0$, it follows from the regularity of $V$ that 
\beaa
 \frac{d}{dt}\tilde\dbE\Big[\big\langle \pa_{x\mu}V(t,X_t,\mathcal{L}_{X_t|\cF^0_t},\tilde X_t)\tilde\d X_t,~\d X_t\big\rangle \Big]\leq 0.
 \eeaa
This is exactly \reff{dI}. 
\qed

\no {\bf Proof of \reff{cdbarI}}. We first apply the It\^o formula \reff{Ito} to obtain 
\begin{equation}\label{eq:ito_double2}
 {d\over dt} {\bar I}(t)= \bar I_1 +\bar I_2 +\bar I_3,
\end{equation}
where, 
\beaa
&&\bar I_1:= \tilde\dbE\bigg[\Big\langle\Big\{\pa_{txx} V(X_t)  + {\widehat\beta^2\over 2} (\tr(\pa_{xx})\pa_{xx} V)(X_t)  +\wh H_p(X_t)^\top \pa_{xxx} V(X_t)\\
&& \qq +\beta^2(\tr(\pa_{x\mu})\pa_{xx}V)(X_t,\tilde X_t)\Big\} \delta X_t, \delta X_t \Big\rangle \bigg],\\
&&\overline{I}_2 :=\bar{\tilde\dbE}\bigg[\Big\langle \Big\{\frac{\beta^2}{2}(\tr(\pa_{\mu\mu})\pa_{xx}V)(X_t,\tilde X_t,\bar X_t)\\
&&\qq +  {\widehat\beta^2\over 2} (\tr(\pa_{\tilde x\mu})\pa_{ xx} V)(X_t,\tilde X_t) +  \wh H_p(\tilde X_t)^\top\pa_{\mu xx} V(X_t,  \tilde X_t)\Big\} \delta  X_t, \delta  X_t\Big\rangle\bigg],\\
&&\overline{I}_3:= 2\bar{\tilde\dbE}\bigg[ \Big\langle\pa_{xx} V(X_t)  \Big\{ \big[\wh H_{px} (X_t) + \wh H_{pp} (X_t)\pa_{xx} V(X_t) \big]\delta  X_t\\
&&\qq \big[ \wh H_{p\rho_1}(X_t, \tilde X_t)+\wh H_{p\rho_2}(X_t,\tilde X_t)\pa_{xx}V(\tilde X_t) + \wh H_{pp}(X_t)\pa_{x\mu} V(X_t, \tilde X_t) \big] \delta \tilde X_t\\
&&\qq + \wh H_{p\rho_2}(X_t,\tilde X_t)\pa_{x\mu}V(\tilde X_t,\bar X_t)\delta \bar X_t\Big\},~ \delta  X_t \Big\rangle \bigg].
\eeaa
On the other hand, applying $\partial_{xx}$ to \reff{master} we obtain  
\begin{equation}\label{paxxVmaster}
0=  (\pa_{xx} \sL V)(t, x, \mu)= \bar J_1 + \bar{J_2},
\end{equation}
\no where 
\begin{align*}
 {\bar J_1}&:= \pa_{txx} V + {\widehat\beta^2\over 2}(\tr(\pa_{xx})\pa_{xx} V) + \wh H_{xx}(x) + 2\wh H_{xp}(x)\pa_{xx} V(x)\\
&+ \pa_{xx} V(x)\wh H_{pp}(x) \pa_{xx} V(x)+\wh H_p(x)^{\top} \pa_{xxx} V(x),\\
 {\bar J_2}&:= \bar{\tilde\dbE}\Big[{\widehat\beta^2\over 2} (\tr(\pa_{\tilde x\mu})\pa_{xx} V)(x,  \tilde  \xi) + \wh H_p( \tilde \xi)^\top\pa_{\mu xx} V(x,  \tilde \xi)\\ 
&+\beta^2(\tr(\pa_{x\mu})\pa_{xx}V)(x,\tilde\xi)+\frac{\beta^2}{2}(\tr(\pa_{\mu\mu})\pa_{xx}V)(x,\bar\xi,\tilde\xi)\Big].
\end{align*}
\normalsize
Evaluate \reff{paxxVmaster} along $(X_t,\mu_t)$, and plug into \reff{eq:ito_double2}, we obtain \reff{cdbarI} straightforwardly. 
\qed

\end{document}